\documentclass[final,3p,times]{elsarticle}
\usepackage{lineno}
\usepackage{graphicx}
\usepackage{amsmath}
\usepackage{natbib}
\usepackage{subcaption}
\usepackage{color}
\usepackage[normalem]{ule m}
\usepackage[shortlabels,inline]{enumitem}
\usepackage{natbib}
\usepackage{tabularx}
\newcolumntype{b}{X}
\newcolumntype{s}{>{\hsize=.5\hsize}X}
\usepackage{bm}
\usepackage{bbm}
\usepackage{multirow}
\usepackage{xcolor}
\usepackage{tikz}
\usepackage[T1]{fontenc}

\colorlet{revision}{purple}

\DeclareRobustCommand\stl{\tikz[baseline]\draw[solid] (0,.5ex)--++(.5,0) ;}
\DeclareRobustCommand\dedash{\tikz[baseline]\draw[densely dashed] (0,.5ex)--++(.5,0) ;}
\DeclareRobustCommand\ddd{\tikz[baseline]\draw[dash dot] (0,.5ex)--++(.5,0) ;}
\DeclareRobustCommand\dd{\tikz[baseline]\draw[dotted] (0,.5ex)--++(.5,0) ;}
\newcommand{\circles}{\protect\tikz{\protect\draw (0.2,0.1) circle (3pt);}}
\newcommand{\squares}{\protect\tikz{\protect\draw (0,0) rectangle (0.2,0.2);}}
\newcommand{\diamonds}{\protect\tikz{\protect\draw[rotate=45] (0,0) rectangle (0.2,0.2);}}
\begin{document}

\begin{frontmatter}

\title{Characterization of the forcing and sub-filter scale terms in the volume-filtering immersed boundary method.}

\author[aff1]{Himanshu Dave}
\author[aff1]{Marcus Herrmann}
\author[aff2]{Peter Brady}
\author[aff1]{M. Houssem Kasbaoui\corref{cor1}}
\cortext[cor1]{Corresponding author, email: houssem.kasbaoui@asu.edu}
\affiliation[aff1]{
             organization={School for Engineering of Matter, Transport and Energy},
             addressline={Arizona State University},
             city={Tempe},
             postcode={85281},
             state={AZ},
             country={USA}}

\affiliation[aff2]{
             organization={CCS-2},
             addressline={Los Alamos National Laboratory},
             city={Los Alamos},
             postcode={87545},
             state={NM},
             country={USA}}

\begin{abstract}
  We present a characterization of the forcing and the sub-filter scale terms produced in the volume-filtering immersed boundary (VF-IB) method by \citet{daveVolumefilteringImmersedBoundary2023}. The process of volume-filtering produces bodyforces in the form of surface integrals to describe the boundary conditions at the interface. Furthermore, the approach also produces unclosed terms called $\tau_\mathrm{sfs}$. The level of contribution from $\tau_\mathrm{sfs}$ on the numerical solution depends on the filter width $\delta_f$. In order to understand these terms better we take take a 2 dimensional, varying coefficient hyperbolic equation shown by \citet{bradyFoundationsHighorderConservative2021}. This case is chosen for two reasons. First, the case involves 2 distinct regions seperated by an interface, making it an ideal case for the VF-IB method. Second, an existing analytical solution allows us to properly investigate the contribution from $\tau_\mathrm{sfs}$ for varying $\delta_f$. The filter width controls how well resolved the interface is. The smaller the filter width, the more resolved the interface will be. A thorough numerical analysis of the method is presented, as well as the effect of $\tau_\mathrm{sfs}$ on the numerical solution. In order to perform a direct comparison, the numerical solution is compared to the filtered analytical solution. Through this we highlight three important points. First, we present a methodical approach to volume filtering a hyperbolic PDE. Second, we show that the VF-IB method exhibits second order convergence with respect to decreasing $\delta_f$ (i.e. making the interface sharper). Finally, we show that $\tau_\mathrm{sfs}$ scales with $\delta_{f}^2$. Large filter widths would require a modeling approach to sufficiently resolve $\tau_\mathrm{sfs}$. However for finer filter widths that have a sufficiently sharp interface, $\tau_\mathrm{sfs}$ can be ignored without any significant reduction in the accuracy of solution.

\end{abstract}


\end{frontmatter}

\section{Introduction}
\label{sec_intro}

Most computational simulations to model physical systems involve bounding surfaces with complex topological interfaces. In the case of fluid flows, examples include airfoils, stirring tanks, turbines and fluidized beds. Accuracy of numerical simulations for such systems hinges on the solver's ability to preserve the fundamental physics, and accurately resolving the interface. Furthermore, considering practical computational costs requires the solver to be robust, scalable and quick. Historically, body-conformal meshes have been the tool of choice to simulate such systems \citep{andersonGridGenerationFlow1994,saadatSemiLagrangianLatticeBoltzmann2020,horneMassivelyparallelUnstructuredOverset2019,koblitzDirectNumericalSimulation2017,vremanStaggeredOversetGrid2017}. While an accurate approach, mesh generation can become a cumbersome process, particularly for complex topologies. In addition, this mesh needs to be regenerated at every timestep for moving interfaces \citep{mavriplisUnstructuredGridTechniques1997}, making it a computationally expensive process in large scale simulations. Immersed boundary methods (IBM), originally proposed by Peskin \citep{laiImmersedBoundaryMethod2000,peskinFluidDynamicsHeart1982,peskinFlowPatternsHeart1972,peskinImmersedBoundaryMethod2002} and later improved by several investigators \citep{uhlmannImmersedBoundaryMethod2005,kempeImprovedImmersedBoundary2012,luoFullScaleSolutionsParticleLaden2007,breugemSecondorderAccurateImmersed2012,glowinskiFictitiousDomainApproach2001,patankarFormulationFastComputations2001} have become an attractive option compared to body-fitted methods for two reasons. First, Cartesian grids are used to represent the fluid, alleviating the need for complex mesh creation around topologically complex interfaces. Second, since the interface is `immersed' within the grid using a cloud of Lagrangian markers, for moving interfaces only the markers need to be transported \citep{uhlmannImmersedBoundaryMethod2005}. The Lagrangian forcing that defines the boundary condition at the interface is transformed onto the Eulerian field using convulations with a regularized Dirac delta \citep{romaAdaptiveVersionImmersed1999}. Despite the popularity of IBMs, certain questions still remain due to the ad-hocness of the forcing term at the immersed boundary, since it does not correspond to any physical term in the Navier-Stokes equations. Furthermore, the development of internal flow within solid bodies and how it affects the hydrodynamic force at the interface is also a point of discussion \citep{kempeImprovedImmersedBoundary2012,uhlmannImmersedBoundaryMethod2005,tschisgaleNoniterativeImmersedBoundary2017}.

Recently \citet{daveVolumefilteringImmersedBoundary2023} presented a novel IB method using the volume filtering technique of \citet{andersonFluidMechanicalDescription1967} called the Volume-Filtered Immersed Boundary (VF-IB) method. Using this approach, they provide sound answers to questions regarding the IB method such as (i) an analytical expression for the immersed boundary forcing term, (ii) Elucidating the role of internal flow, and (iii) a more accurate approach to calculating the Lagrangian marker volumes. The transport equations are derived by filtering the original point-wise equations in order to obtain a new set of filtered equations. The boundary conditions which are normally imposed on the fluid-solid interface are converted into bodyforces that apply on the right-hand side of the filtered transport equations as surface integrals. The process is mathematically and physically rigorous, and does not depend on any numerical considerations. They show that by doing this they remove any ad-hoc numerical fixes such as retraction of the immersed boundary to get accurate hydrodynamic forces \citep{breugemSecondorderAccurateImmersed2012}. Lastly, they show that the Lagrangian marker volumes depend on the local topology of the interface, the choice of filter kernel and the local curvature of the interface. Accurately calculating the Lagrangian marker volume can help improve the calculation of the hydrodynamic force due to more accurate interpolation and extrapolation operations. They also show an efficient procedure to compute the volume fractions of the different regions using a Poisson equation, required to accurately calculate the stresses at the interface.

While the VF-IB method helped answer several longstanding questions, it also has some questions that arise out of the volume-filtering procedure which require careful consideration. The process of volume-filtering the point-wise equations produces unclosed terms, including a subfilter scale tensor $\tau_\mathrm{sfs}$, similar in spirit to the Large Eddy Simulations \citep{daveVolumefilteringImmersedBoundary2023}. The level of contribution from $\tau_\mathrm{sfs}$ depends on the size of the filter width. \citet{daveVolumefilteringImmersedBoundary2023} ignore this term and show good accuracy for sharp interface resolution. However, this value cannot be ignored as the filter width is increased. In order to make the VF-IB method scalable, simulations need to be accurately run at coarse resolution, thus making it important to understand the role of $\tau_\mathrm{sfs}$ in comparison to the other terms in the filtered equation.

To understand $\tau_\mathrm{sfs}$ and $F_\mathrm{IB}$, we inspect the VF-IB method using an equation that has an analytical solution. This can help us properly understand the contribution of $\tau_\mathrm{sfs}$ and other terms for varying parameters, allowing us to properly characterize the method. To do this we use a varying coefficient hyperbolic equation \citep{bradyFoundationsHighorderConservative2021}. This is an archetypal problem to solve using the VF-IB method and ideal to help characterize it for two reasons. First, the problem at hand has two distinct regions seperated by a circular interface, making it a perfect candidate for the VF-IB method. Second, an existing analytical solution allows us to properly perform an error analysis and understand the effects of $\tau_\mathrm{sfs}$ and other terms in the filtered equation for varying parameters.

In section \ref{sec:sec_test_case}, we introduce the test case, show the computational domain, and investigate the analytical solution. In section \ref{sec:sec_vf} we introduce the volume filtering procedure and the prerequisites of the method. With the mathematical framework in place, we show the derivation of the filtered equations from the point wise equations in section \ref{sssec:derivation_vfib}. In section \ref{sssec:vf_compute}, we show how the volume fraction computation is performed. Accurately calculating the volume fraction is neccessary to correctly resolve the interface and help distinguish between the two regions. The temporal and spatial schemes used in the numerical method are shown in section \ref{sec:numerical}. We then conduct a thorough numerical analysis to look at the the order of convergence of the VF-IB method. The results are addressed in section \ref{sec:sec_results}. Section \ref{sssec:ap_analysis} shows an Apriori analysis of the case. We examine the filtered analytical solution which is the basis of comparison to the numerical solution produced by the VF-IB method in order to perform a more direct comparison. Furthermore, we also examine the different terms within the filtered equation and how they compare with each other. This section also explores how the magnitude of $\tau_\mathrm{sfs}$ varies with varying filter width. Additionally, we also show how the subgrid resolution $\delta_f/\Delta x_f$ affects the accuracy of the filtered quantity. We then conduct an aposteriori analysis in section \ref{sssec:apos_analysis}. Here we compare the numerical solution produced by the VF-IB method and the filtered analytical solution, both for varying filter width and varying grid resolution while keeping $\tau_\mathrm{sfs}$ turned off. Lastly we investigate how including $\tau_\mathrm{sfs}$ affects the solution for varying filter width. Finally, we give the concluding remarks in section \ref{sec:sec_conclusion}.

\section{Varying coefficient hyperbolic equation}
\label{sec:sec_test_case}
\begin{figure}
	\centering
	\includegraphics[width=0.7\linewidth]{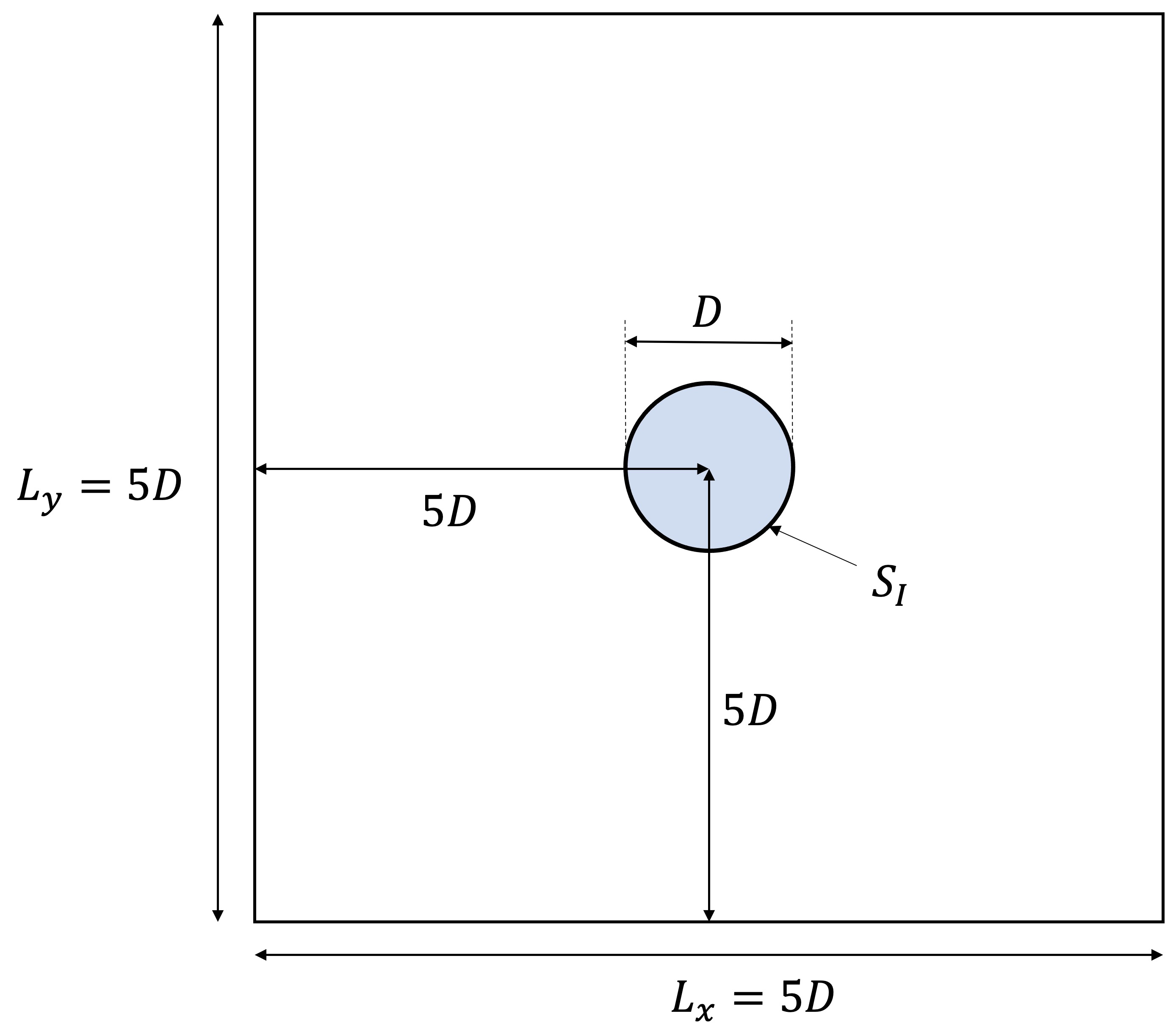}
	\caption{Computational domain for the test case of a varying coefficient hyperbolic equation similar to the domain used by \citet{bradyFoundationsHighorderConservative2021}.}
	\label{fig:schematic_testcase}
\end{figure}
In order to test the VF-IB method, we consider the two-dimensional test case of a varying coefficient hyperbolic equation shown in \citet{bradyFoundationsHighorderConservative2021}. The governing equation in non-dimensional form is given as follows,
\begin{eqnarray}
  \frac{\partial u}{\partial t} + \nabla G\cdot\nabla u &=& 0,\label{eq:GE_1}
\end{eqnarray}
where $G(x,y)$ is,
\begin{eqnarray}
  G(x,y)&=&\sqrt{\left(x-x_c\right)^2+\left(y-y_c\right)^2}-r,\label{eq:GE_2a}\\
  \Omega_2 &=& \{(x,y)\mid G(x,y)<0\},\label{eq:GE_2b}\\
  S_I &=& \{(x,y)\mid G(x,y)=0\}.\label{eq:GE_2c}
\end{eqnarray}
In equation (\ref{eq:GE_2a}), $x_c$ and $y_c$ are the center coordinates of the circle and $r$ is the radius. $S_I$ represents the interface between the two regions and $\Omega_2$ represents the region inside the circle. The computational domain is shown in figure \ref{fig:schematic_testcase}. we take a domain size of $L_x = L_y = 2$ and $D = 0.4$. The circle is located at the center of the computational domain. The initial and boundary conditions are given by,
\begin{eqnarray}
  u(x,y,t=0) &=& \mathrm{sin}\left(2\pi G\right),\label{eq:GE_3a}\\
  u_I(x,y,t)\mid_{G=0} &=& -\mathrm{sin}\left(2\pi t\right),\label{eq:GE_3b}
\end{eqnarray}
and the boundary conditions at the edge of the computational domain is an outflow condition. The analytical solution for the case is a circular pulse radiating out from the circle with a period of 1 such that,
\begin{eqnarray}
  u(x,y,t) = \mathrm{sin}\left(2\pi\left(G-t\right)\right).\label{eq:GE_4}
\end{eqnarray}

\begin{figure}
	\centering
	\includegraphics[width=0.9\linewidth]{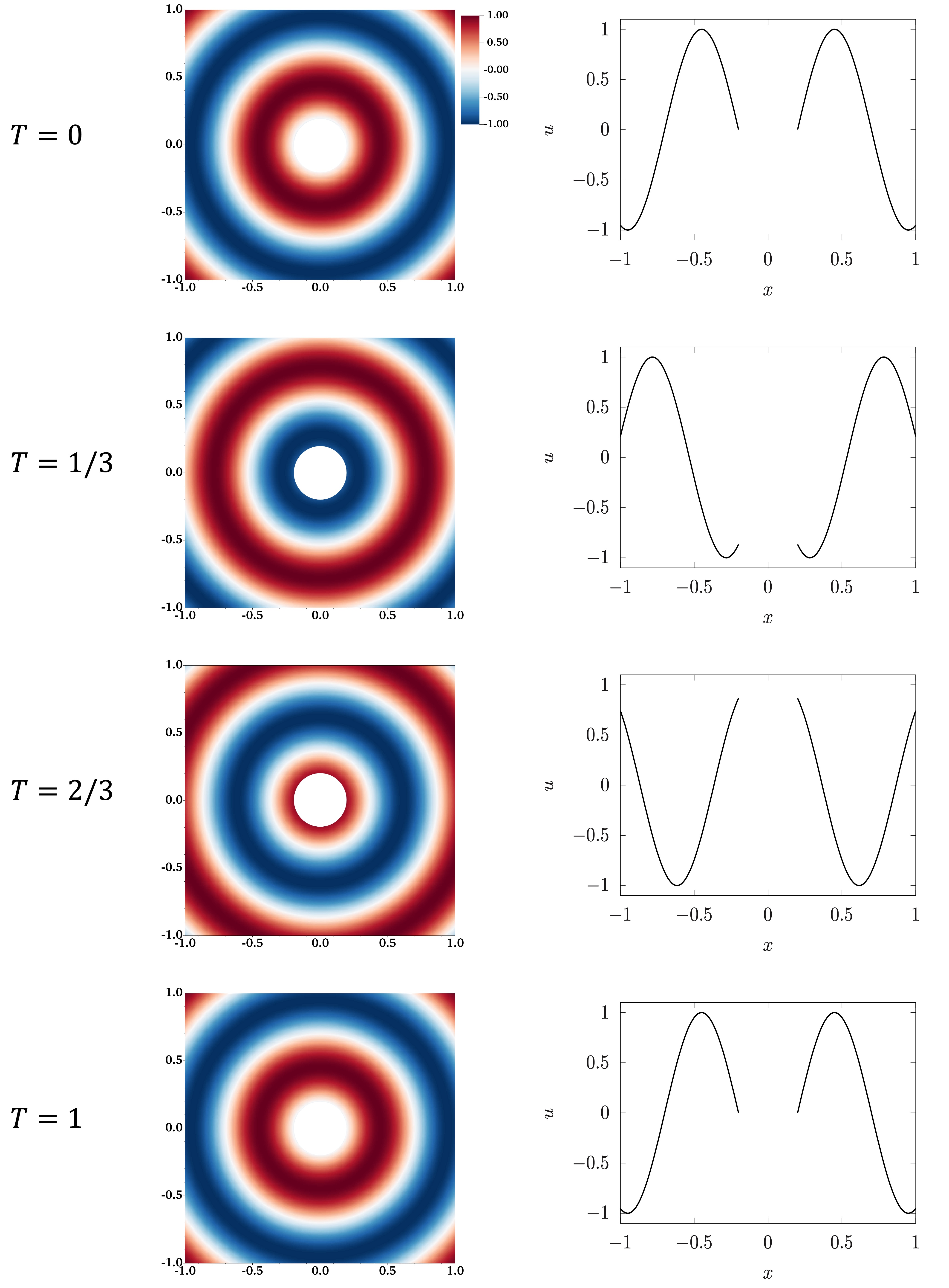}
	\caption{Analytical solution for a varying coefficient hyperbolic equation. The isocontours of $u$ at different time periods $T$ (left) and the value of $u$ vs $x$ in the horizontal direction, centered in the vertical direction (right). There is no velocity field within the region inside the circle.}
	\label{fig:analytical}
\end{figure}

Figure \ref{fig:analytical} shows the analytical solution for one period of the circular radiating pulse. A clear circular pulse radiating out from the circle is shown. Furthermore, we display the graph of $u$ vs $x$ along the horizontal direction centered in the vertical direction. The radiating pulse is symmetric in all directions.

\section{Volume-filtering approach}
\label{sec:sec_vf}
\begin{figure}
	\centering
	\includegraphics[width=0.9\linewidth]{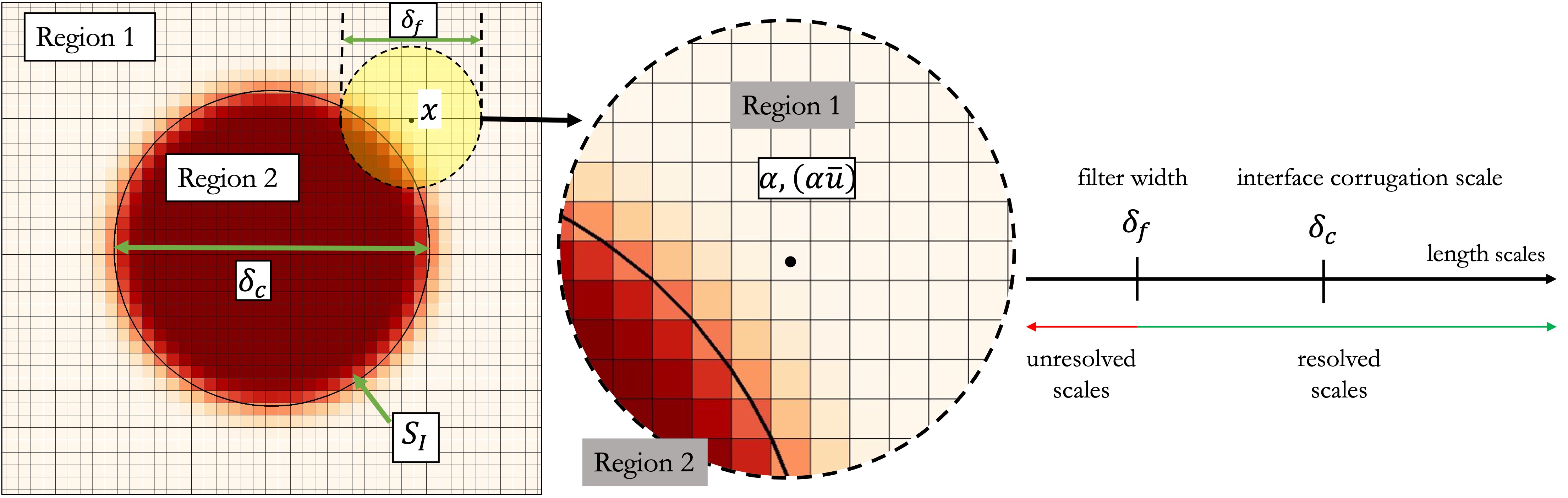}
	\caption{ Illustration of the volume-filtering approach. Filtering the point-wise fields allows the extraction of the volume fraction $\alpha$, and the averaged point-wise fields $(\alpha\overline{u})$ for region 1. The averaging is performed within a region of size $\delta_f$. The immersed boundary is well resolved when the characteristic corrugation scale $\delta_c$ of the interface is much larger than the filter width $\delta_f$. ($\delta_f$ is not to scale in the figure, but rather shown much larger for easier understanding of the concept of volume-filtering).}
	\label{fig:schematic_VF}
\end{figure}

In this section, we show the volume filtering approach by \citet{andersonFluidMechanicalDescription1967} and how it can be applied in conjunction with the immersed boundary method. Figure \ref{fig:schematic_VF} shows an illustration of the volume filtering approach. Here two distinct regions are seperated by an immersed boundary surface $S_I$. While there is a sharp discontinuity between the regions analytically, their respective volume fraction fields are a smeared indicator function. The degree of smearing is dependent on the size of the filter width $\delta_f$. For illustration, $S_I$ in figure \ref{fig:schematic_VF} is chosen by the iso-level $\alpha = 0.5$, where $\alpha$ represents the volume fraction related to region 1.

To understand how the volume filtering process works, let us take the point-wise quantity $u$ defined in region 1. We now take a point $\bm{x}$ close to the immersed boundary surface but within region 1 as shown in figure \ref{fig:schematic_VF}. $u$ can then be filtered at point $\bm{x}$. The size of the area within which the filtering process takes place to obtain an averaged quantity at $\bm{x}$ depends on the size of the filter width $\delta_f$. Through this process, we obtain the volume fraction $\alpha$, and the averaged point-wise quantity $(\alpha\overline{u})$ at $\bm{x}$. While $u$ only exists in region 1, the volume-filtered field $(\alpha\overline{u})$ exists everywhere. As we move past the interface within region 2, $(\alpha\overline{u})$ smoothly decays to zero. There exists a region $\delta_f/2$ away from the interface within region 2 where $(\alpha\overline{u})$ is still non-zero. This is because for points within region 2 close to the interface, there exists a partial area under the filter width that is within region 1. In the case of figure \ref{fig:schematic_VF}, we show the volume fraction field $\alpha$. The size of the region considered when averaging a point-wise quantity depends on the filter width $\delta_f$. Hence, smaller the filter width, the faster $\alpha$ and $(\alpha\overline{u})$ will decay to zero as we move into region 2. Furthermore, how well resolved the immersed boundary is dependent on the ratio between the interface corrugation length scale $\delta_c$ and the filter kernel width $\delta_f$. In order to have a well resolved immersed boundary we require that $\delta_f \ll \delta_c$.

In order to formalize the idea of volume filtering, we consider a filter kernel $g$ that satisfies,
\begin{align}
    \iiint_{\bm{y}\in \mathbbm{R}^3} g(\bm{y})dV&\;=\;1,&\text{(unitary)}\label{eq:filter_def}\\
    g(-\bm{y})&\;=\;g(\bm{y}),&\text{(symmetric)}\\
    g(\bm{y})&\;=\;0\ \mbox{if}\ ||\bm{y}||\geq \delta_f/2.&\text{(compact)}
\end{align}
The property of symmetry is important since it helps eliminate artificial anisotropy. Compactness helps for fast numerical integration of $g$ on surfaces. The integration is considered over the entire space.

The volume fraction at any arbitrary location $\bm{x}$ is given by
\begin{align}
  \alpha(\bm{x},t)&\;=\;\iiint_{\bm{y}\in \mathbbm{R}^3}\mathbbm{1}(\bm{y},t)g(\bm{x}-\bm{y})dV\label{eq:RT1_2}.
\end{align}
Here $\mathbbm{1}(\bm{y},t)$ is an indicator function equal to 1 if $\bm{y}$ is in region 1 and 0 otherwise. $\alpha(\bm{x})$ represents the total region 1 that exists within the support of the filter kernel. If $\bm{x}$ is far away from region 2 such that the entire area under the filter kernel support is exclusively within region 1, the value of $\alpha(\bm{x}) = 1$. If the probe is moved closer to the interface, such that the area within the support of the kernel is under both the regions, the value of $\alpha(\bm{x})$ will be somewhere between $0<\alpha(\bm{x})<1$. Through this, we eliminate the discontinous effects across the interface by smoothing out the volume fraction field based on the scale of the filter width $\delta_f$ chosen. In a similar approach we obtain the filtered quantity for the point-wise quantity $u$ that exists at any arbitrary point $\bm{x}$, which is defined as,
\begin{eqnarray}
  \alpha(\bm{x},t)\overline{u}(\bm{x},t)&=&{ \iiint_{\bm{y}\in\mathbbm{R}^3} \mathbbm{1}(\bm{y},t)u(\bm{y},t)g(\bm{x}-\bm{y})dV}.
  \label{eq:filtered_lambda}
\end{eqnarray}
The volume-filtered quantity $(\alpha\overline{u})$ is continous and exists everywhere. This value tends smoothly to 0 a distance $\delta_f/2$ away from the interface within region 2.

Following the work of \citet{andersonFluidMechanicalDescription1967} we obtain the volume-filtered equations from the governing point-wise equations. In order to do this, we apply the filtering operations to the point-wise equations. Using the property of symmetry and the divergence theorem, we obtain the filtered gradient and time derivative operators,
\begin{eqnarray}
  \alpha(\bm{x}) \overline{\nabla u}(\bm{x}) &=&\nabla (\alpha\overline{u}) - \iint_{\bm{y}\in S_I} \bm{n}u(\bm{y},t) g(\bm{x}-\bm{y})dS,\label{eq:identity1}\\
  \alpha(\bm{x}) \overline{\frac{\partial u}{\partial t}}(\bm{x}) &=&\frac{\partial (\alpha \overline{u})}{\partial t} + \iint_{\bm{y}\in S_I} (\bm{n}\cdot\bm{u}_{\mathrm{IB}})u(\bm{y},t) g(\bm{x}-\bm{y})dS.\label{eq:identity3}
\end{eqnarray}
Here $\bm{n}$ is the normal vector pointing from region 2 to region 1 at the interface $S_I$. In equation (\ref{eq:identity3}), $\bm{u}_\mathrm{IB}$ is the velocity of the immersed boundary, not to be confused with the boundary condition $u_I$, at the interface. For static boundaries, this value will be zero. The process of volume filtering removes the notion of a boundary, since the filtered quantities exist everywhere in $\mathbbm{R}^3$. The information of the boundary conditions emerge in the surface integrals that arise in equations (\ref{eq:identity1}) and (\ref{eq:identity3}).

\subsection{Derivation of the volume-filtered equations}
\label{sssec:derivation_vfib}

We now take equation (\ref{eq:GE_1}) and perform volume-filtering operations on it to obtain the filtered governing equation. Filtering the first term in equation (\ref{eq:GE_1}) leads to,
\begin{eqnarray}
	\alpha\overline{\frac{\partial u}{\partial t}}&=&\frac{\partial}{\partial t}\left(\alpha\overline{u}\right).\label{eq:VF_1}
\end{eqnarray}
For a fixed interface, volume-filtering the time derivative leads to no surface integrals since the interface velocity, $\bm{u}_\mathrm{IB}$ is zero. In order to volume filter the second term in equation (\ref{eq:GE_1}), we expand the filtered quantity as follows,
\begin{eqnarray}
	\alpha\overline{\nabla G\cdot\nabla u}&=&\alpha\overline{\nabla G}\cdot\overline{\nabla u} + \left(\alpha\overline{\nabla G\cdot\nabla u} - \alpha\overline{\nabla G}\cdot\overline{\nabla u}\right).\label{eq:VF_2}
\end{eqnarray}
Volume filtering the first term on the right hand side of equation (\ref{eq:VF_2}) leads to,
\begin{eqnarray}
	\alpha\overline{\nabla G}\cdot\overline{\nabla u}&=&\overline{\nabla G}\cdot\nabla\left(\alpha\overline{u}\right) - \overline{\nabla G}\cdot\iint_{\bm{y}\in S_I}u_I\bm{n}g(\bm{x}-\bm{y})dS.\label{eq:VF_3}
\end{eqnarray}
Here $u_I$ is the prescribed quantity at the interface. In equation (\ref{eq:VF_2}), we see the emergence of the sub-filter scale term in the last term on the right hand side such that,
\begin{eqnarray}
	\tau_{\mathrm{sfs}}&=&\alpha\overline{\nabla G\cdot\nabla u} - \alpha\overline{\nabla G}\cdot\overline{\nabla u}.\label{eq:VF_4}
\end{eqnarray}

This term is similar in fashion to the sub-grid scale terms that emerge from the Large-Eddy Simulation (LES) equation of the Navier-Stokes equation. Combining all the terms, we get the final volume-filtered governing equation for a varying coefficient hyperbolic equation as,
\begin{eqnarray}
	\frac{\partial}{\partial t}\left(\alpha\overline{u}\right) + \overline{\nabla G}\cdot\nabla\left(\alpha\overline{u}\right)&=& F_I - \tau_{\mathrm{sfs}}.\label{eq:VF_5}
\end{eqnarray}
Here $F_I = \overline{\nabla G}\cdot\iint_{\bm{y}\in S_I}u_I\bm{n}g(\bm{x}-\bm{y})dS$. The boundary conditions at the interface arise as surface integrals.

\subsection{Volume-fraction computation}
\label{sssec:vf_compute}

Computing the volume fraction of the different regions within the computational domain serves two goals pertaining to the test case presented above: (i) distinguishing between the interior and exterior points, (ii) computing the total volume occupied by the immersed solid. Additionally, the process of volume filtering produces the volume filtered quantity in the form of $(\alpha\overline{u})$. In certain cases we may be required to extract $\overline{u}$. This requires us to compute the volume fraction in order to perform the operation $(\alpha\overline{u})/\alpha$.

We first compute the region 1 indicator function within the field such that,
\begin{eqnarray}
  \mathbbm{1}(\bm{x},t)=\begin{cases}
    1, & \text{if $\sqrt{(x-x_c)^2+(y-y_c)^2}>R$}.\\
    0, & \text{otherwise}.
  \end{cases}
\end{eqnarray}
where $R$ is the radius of the circle. The center of the circle is defined by $x_c$ and $y_c$. The indicator function is then filtered, in order to obtain the volume fraction such that,
\begin{align}
  \alpha(\bm{x},t)&\;=\;\iiint_{\bm{y}\in \mathbbm{R}^3}\mathbbm{1}(\bm{y},t)g(\bm{x}-\bm{y})dV\label{eq:vfp_compute}.
\end{align}

\begin{figure}\centering
  \includegraphics[width=\linewidth]{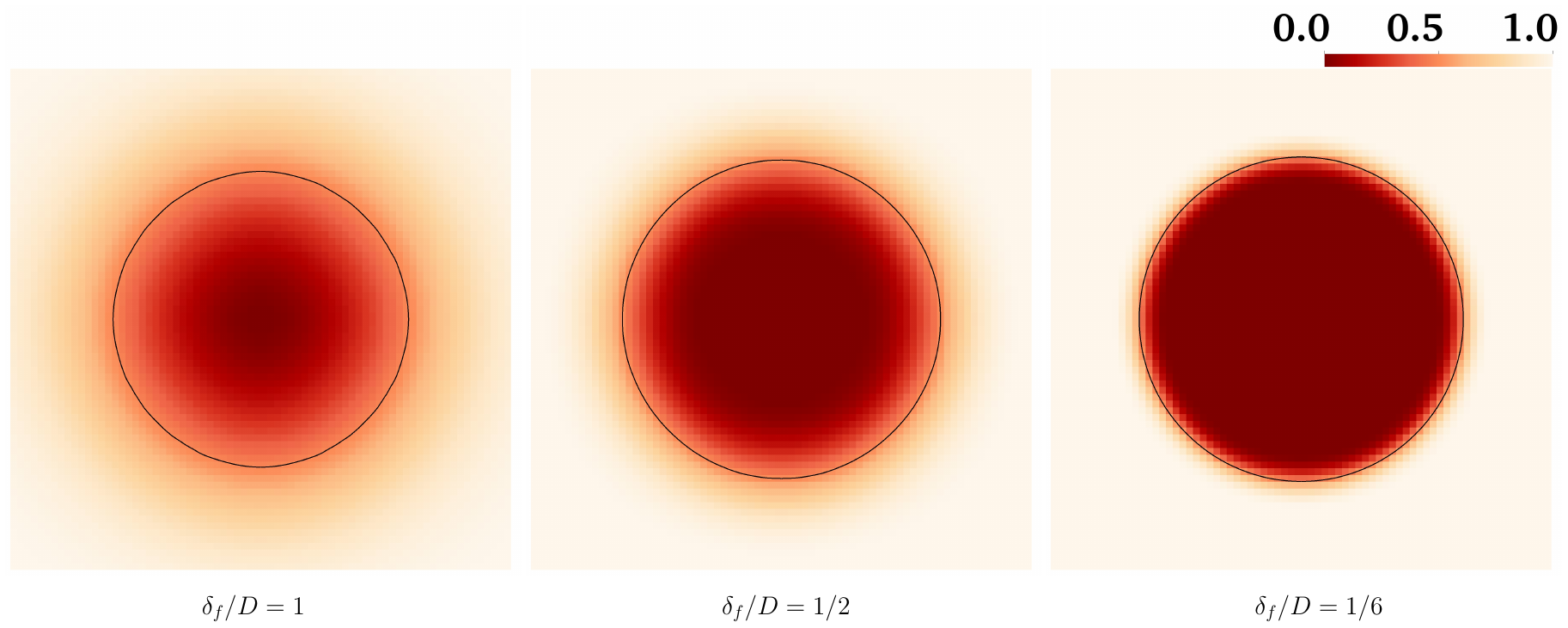}
	\caption{Volume fraction $\alpha$, at three different filter widths with respect to the circle diameter ($\delta_f/D = 1$, $1/2$ and $1/6$). The black contour line represents the Immersed Boundary (IB) surface located at $\alpha = 0.5$ in the limit $\delta_f/D \rightarrow 0$.\label{fig:vfp}}
\end{figure}

Figure \ref{fig:vfp} shows the volume fraction field at 3 different $\delta_f/D$ values of $1$, $1/2$ and $1/6$. Here, we show the diffuse nature of the volume fraction field. The distance it takes to go from $\alpha = 1$ to $\alpha = 0$ depends on the filter kernel support. The larger the filter width, as in the case of $\delta_f/D = 1$ shown in figure \ref{fig:vfp}, the greater the smearing. This is visually observed when looking at the transition range from $\alpha = 1\rightarrow 0$ for the $\delta_f/D = 1$ case. As we make the filter width smaller, i.e. $\delta_f/D = 1/6$, the transition occurs within a shorter distance. Hence, the sharpness of the filtered quantities depends on the filter width chosen. The smaller the filter width with respect to the circle, the sharper the representation of the interface. Furthermore, we also observe that as $\alpha$ transitions from one to zero, there does exist a region $\delta_f/2$ away from the interface but within region 2 where $\alpha$ is non-zero. In order to highlight the smearing effect of volume filtering more clearly we show a contour line at $\alpha = 0.5$ which is chosen to represent the IB surface in the limit $\delta_f/D \rightarrow 0$.

This approach works well for the case chosen due to the simplicity of the geometry together with having a fixed interface. For more complex topological surfaces and/or moving interfaces a more sophisticated approach should be chosen to compute the volume fraction. \citet{daveVolumefilteringImmersedBoundary2023} presents an approach that entails solving a Poisson equation in order to compute the volume fraction such that,
\begin{eqnarray}
	\nabla^{2}\alpha(\bm{x},t) = \nabla\cdot\iint_{\bm{y}\in S_I}\bm{n}g(\bm{x}-\bm{y})dS.\label{eq:VF_10}
\end{eqnarray}
This is similar to the computation of the phase-indicator functions with the front-tracking method of \citet{unverdiFrontTrackingMethodViscous1992}. \citet{daveVolumefilteringImmersedBoundary2023} choose an algebraic multigrid solver with Dirichlet boundary conditions to solve equation (\ref{eq:VF_10}). The volume fraction is computed everywhere in the domain for simplicity. However, it would suffice to solve the Poisson equation in a narrow band with a thickness that is equal to the filter width $\delta_f$ similar to what is done by \citet{unverdiFrontTrackingMethodViscous1992}. In the case of a fixed interface, the Poisson solver would only need to be solved once.

\section{Numerical implementation}
\label{sec:numerical}

The VF-IB method is implemented in a library called LEAP. In this section we describe the numerical implementation of the VF-IB approach in LEAP pertaining to the test case shown in section \ref{sec:sec_test_case}.

\subsection{Spatial discretization of the interface}
\label{sssec:spatial}
We will now examine equation (\ref{eq:VF_5}) and focus on the discretization of the interface. The IB interface is generated using a mesh of `m' discrete elements with a surface area $A_m$, having a centroid $\bm{x}_m$ and an outward pointing normal $\bm{n}_m$. The IB forcing can then be written as a sum of all discrete contributions from each element of the mesh.

The forcing on the interface is performed by the surface integral term in equation (\ref{eq:VF_5}),
\begin{eqnarray}
  F_{I}(\bm{x})&=&\overline{\nabla G}(\bm{x})\cdot\sum_{m=1}^{N}\iint_{\bm{y}\in S_I}u_I\bm{n}\mid_{\bm{y}}g(x-y)dS.\label{eq:SD_1}
\end{eqnarray}
Using the midpoint rule and assuming a typical mesh width of $O(\Delta x)$, equations (\ref{eq:SD_1}) leads to,
\begin{eqnarray}
  F_{I}(\bm{x})&=&\overline{\nabla G}(\bm{x})\cdot\sum_{m=1}^{N}\{\left(u_I\bm{n}\right)\mid_{\bm{x}_m}g(\bm{x}-\bm{x}_m)A_m\}.\label{eq:SD_3}
\end{eqnarray}
Here $\bm{x}_m$ is the centroid location of the mesh element $S_m$ and its surface area is $A_m$. The centroids are portrayed as Lagrangian forcing points as shown in \citet{daveVolumefilteringImmersedBoundary2023} and many other IB methods. The interface forcing is then extrapolated onto the Eulerian field. Details of the extrapolation procedure are shown in \citet{daveVolumefilteringImmersedBoundary2023}.

\subsection{Temporal discretization}
\label{sssec:temporal}

The time integration scheme is based on a strong-stability preserving Runga-Kutta 3rd order (SSP-RK3) scheme. The steps below describe the update from $t^n$ to $t^{n+1}$.

\textbf{Step 1:} At this step the time integration loop is started. The first step is to make $u^{(0)}$ equal to $u^{n}$,
\begin{eqnarray}
  (\alpha\overline{u})^{(0)} = (\alpha\overline{u})^{n}.\label{T_1}
\end{eqnarray}

\textbf{Step 2:} Next, we compute the immersed boundary forcing term, i.e. equation (\ref{eq:SD_3}), for step 1 of the SSP-RK3 scheme,
\begin{eqnarray}
  F^{n}_{I}(\bm{x})&=&\overline{\nabla G}\cdot\sum_{m=1}^{N}\left\{\left(u^{n}_{I,m}\bm{n}\right)\mid_{\bm{x}_m}g(\bm{x}-\bm{x}_m)A_m\right\},\label{eq:T_2}
\end{eqnarray}
The updated immersed boundary forcing is then extrapolated onto the Eulerian field in order to proceed with the time integration scheme. The first step velocity is computed as,
\begin{eqnarray}
  (\alpha\overline{u})^{(1)} = \Delta t\left(-\overline{\nabla G}\cdot\nabla(\alpha\overline{u})^{(0)} + F^{n}_{I} - \tau_{\mathrm{sfs}}^n\right) + (\alpha\overline{u})^{(0)}.\label{eq:T_4}
\end{eqnarray}

\textbf{Step 3:} In similar fashion to expressions (\ref{eq:T_2}), we compute the immersed boundary forcing term for step 2 of the SSP-RK3 scheme for time $t^{n+1}$. The second step velocity is computed as,
\begin{eqnarray}
  (\alpha\overline{u})^{(2)} &=& \frac{3}{4}(\alpha\overline{u})^{(0)} + \frac{1}{4}\left\{\Delta t\left(-\overline{\nabla G}\cdot\nabla\left(\alpha\overline{u}\right)^{(1)} + F^{n+1}_{I} - \tau_{\mathrm{sfs}}^{n+1}\right) + (\alpha\overline{u})^{(1)}\right\}.\label{eq:T_5}
\end{eqnarray}

\textbf{Step 4:} Finally, we compute the final velocity $u^{n+1}$ as,
\begin{eqnarray}
  (\alpha\overline{u})^{n+1} &=& \frac{1}{3}(\alpha\overline{u})^{(0)} + \frac{2}{3}\left\{\Delta t\left(-\overline{\nabla G}\cdot\nabla\left(\alpha\overline{u}\right)^{(2)} + F^{n+1/2}_{I} - \tau_{\mathrm{sfs}}^{n+1/2}\right) + (\alpha\overline{u})^{(2)}\right\}.\label{eq:T_6}
\end{eqnarray}

\section{Results}
\label{sec:sec_results}

In this section we show the comparison of the results obtained from the filtered analytical solution and solution obtained using the VF-IB method. The domain and cylinder radius are kept consistent with what is shown in section \ref{sec:sec_test_case}. First we conduct an apriori analysis (section \ref{sssec:ap_analysis}) to understand the difference between the analytical solution shown in section \ref{sec:sec_test_case} and the filtered analytical solution. We then examine the various terms within the filtered governing equation to show how each term plays a role. Furthermore, we show how $\tau_\mathrm{sfs}$ changes as $\delta_f/D$ is varied. Secondly in section \ref{sssec:apos_analysis}, we perform an error analysis and obtain and order of convergence for the VF-IB method for both spatial quantities, $\delta_f/\Delta x$ and $\delta_f/D$. We also investigate how including $\tau_\mathrm{sfs}$ in the solution affects it for varying filter width. In order to distinguish between the different quantities, $u_e$ will be the unfiltered analytical solution, $(\alpha\overline{u})_e$ is the filtered analytical solution and $(\alpha\overline{u})$ is the numerical solution produced using the VF-IB method.

\subsection{Apriori analysis}
\label{sssec:ap_analysis}
\begin{figure}
	\centering
	\includegraphics[width=\linewidth]{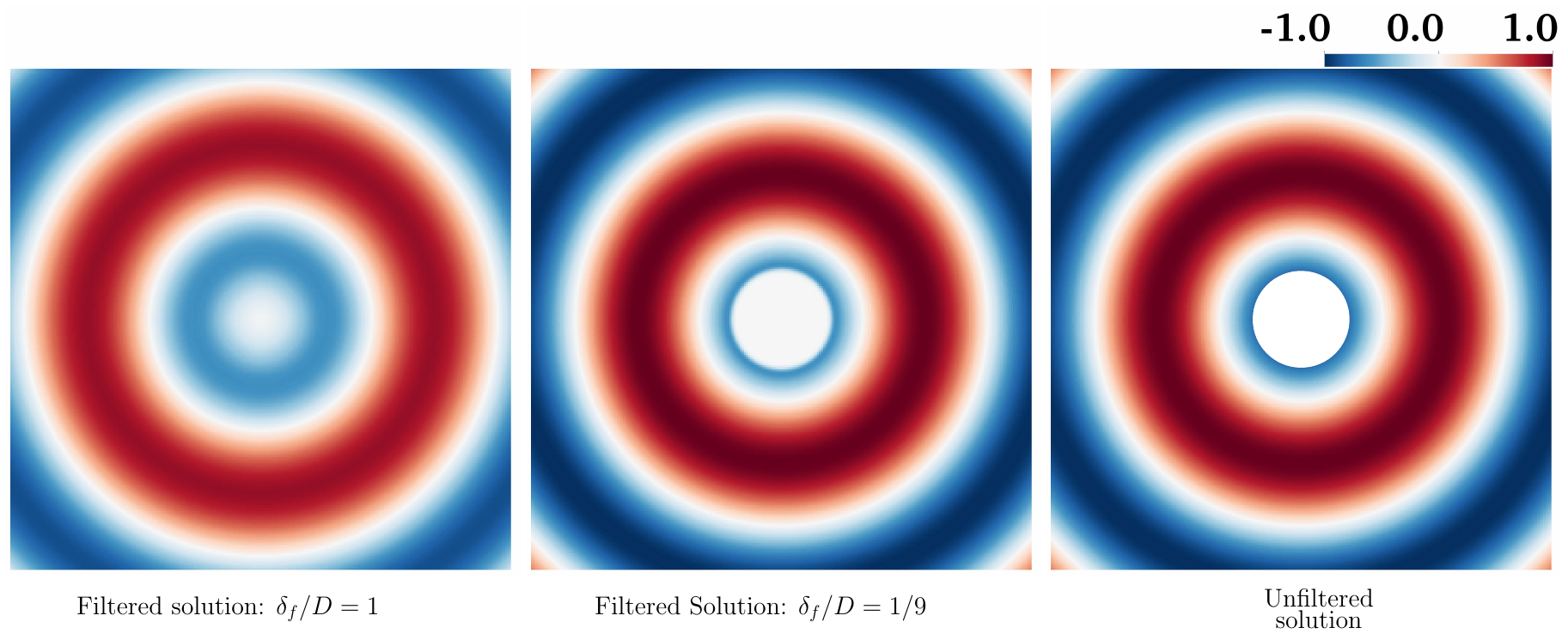}
	\caption{Isocontours of the filtered solution $(\alpha\overline{u})_e$, and the unfiltered solution $u_\mathrm{e}$ at time period $T = 1/4$. At this time the forcing is at its absolute maximum value. $(\alpha\overline{u})_e$ is shown at 2 different filter resolution of $\delta_f/D = 1$ and $\delta_f/D = 1/9$.}
	\label{fig:iso_ufexact}
\end{figure}
The accuracy of the filtered quantity in comparison to the unfiltered value is dependent on the filter width $\delta_f$. The smaller the filter width, the sharper the representation of the filtered quantity is and the more closer it is to the unfiltered value. As the filter width is increased, the filtered field becomes more smeared and digresses from the unfiltered value. In order to show this, figure \ref{fig:iso_ufexact} shows the isocontours for $(\alpha\overline{u})_e$ at two different filter resolutions of $\delta_f/D = 1$ and $\delta_f/D = 1/9$. The results are at the same time, $T = 1/4$. We also show the unfiltered solution $u_e$ at the same time for comparison. From the instantaneous visualizations, it is immediately clear that $(\alpha\overline{u})_e$ exists everywhere in space. The quantity smoothly decays to zero as we move past the interface towards the center of the circle. The distance past the interface at which $(\alpha\overline{u})_e$ is zero is $\delta_f/2$. In comparison, $u_e$ only exists in region 1 and does not go past the interface. As shown in figure \ref{fig:iso_ufexact}, at a filter width of $\delta_f/D = 1$ the smearing of the field is greater than at $\delta_f/D = 1/9$. Hence, for larger filter widths the transition region will be larger, as we move from region 1 into region 2. The smaller the filter width the sharper the field looks and the closer it is to the unfiltered solution.

\begin{figure}\centering
  \begin{subfigure}{0.45\linewidth}
    \includegraphics[width=\linewidth]{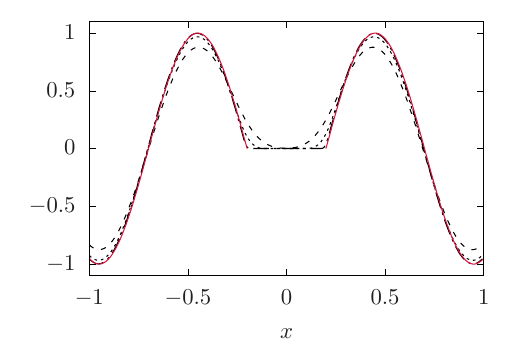}
    \caption{$T = 0$\label{fig:fig6a}}
  \end{subfigure}
  \begin{subfigure}{0.45\linewidth}
    \includegraphics[width=\linewidth]{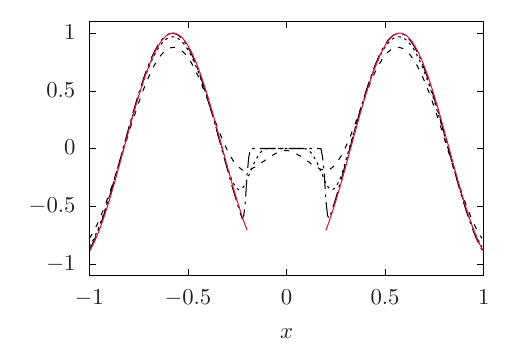}
    \caption{$T = 1/8$\label{fig:fig6b}}
  \end{subfigure}
	\begin{subfigure}{0.45\linewidth}
    \includegraphics[width=\linewidth]{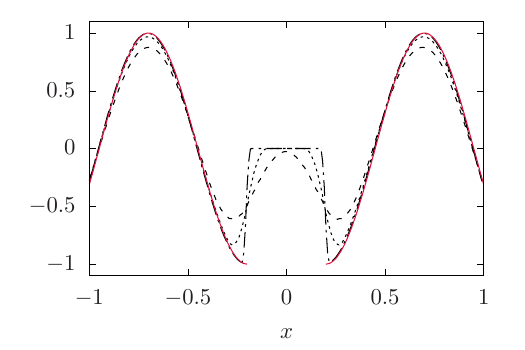}
    \caption{$T = 1/4$\label{fig:fig6c}}
  \end{subfigure}
	\caption{Value of $(\alpha\overline{u})_e$ vs $x$ at $T = 0$, $T = 1/8$ and $T = 1/4$. The result shown is a line cut horizontally and centered in the vertical axis. We show three different filter resolutions of $\delta_f/D = 1$ (\dedash), $\delta_f/D = 1/2$ (\dd) and $\delta_f/D = 1/9$ (\ddd). We also show $u_\mathrm{e}$ (\textcolor{red}{\stl}) for reference.\label{fig:figure6}}
\end{figure}

In order to quantitatively highlight the differences between the filtered and unfiltered solution, figure \ref{fig:figure6} shows $(\alpha\overline{u})_e$ at different filter resolutions alongside the unfiltered solution $u_\mathrm{e}$. We choose three different filter widths ranging from $\delta_f/D = 1$ to $\delta_f/D = 1/9$. $(\alpha\overline{u})_e$ and $u_\mathrm{e}$ are compared along the horizontal axis, centered in the vertical axis. We show the graph at three instances in time. Figure \ref{fig:fig6a} is at $T = 0$, corresponding to the time when the forcing is at its absolute minimum value. In figure \ref{fig:fig6b}, we show the results at $T = 1/8$, when the forcing is half of its maximum absolute value. Finally, we show the results when the forcing term is at its absolute maximum value at $T = 1/4$ in figure \ref{fig:fig6c}. Again, we see that while the filtered quantity exists past the boundary, the unfiltered solution only exists outside the circle. Furthermore, we notice that as we reduce the filter width, the value of $(\alpha\overline{u})_e$ becomes sharper. The results obtained using a filter width of $\delta_f/D = 1/9$ follows the unfiltered solution well. However, with $\delta_f/D = 1$, we see greater smearing of the results. This is particularly visible close to the boundary ($x = -0.2$ and $x = 0.2$) and near the peaks of the solution.

\begin{figure}
	\centering
	\includegraphics[width=0.9\linewidth]{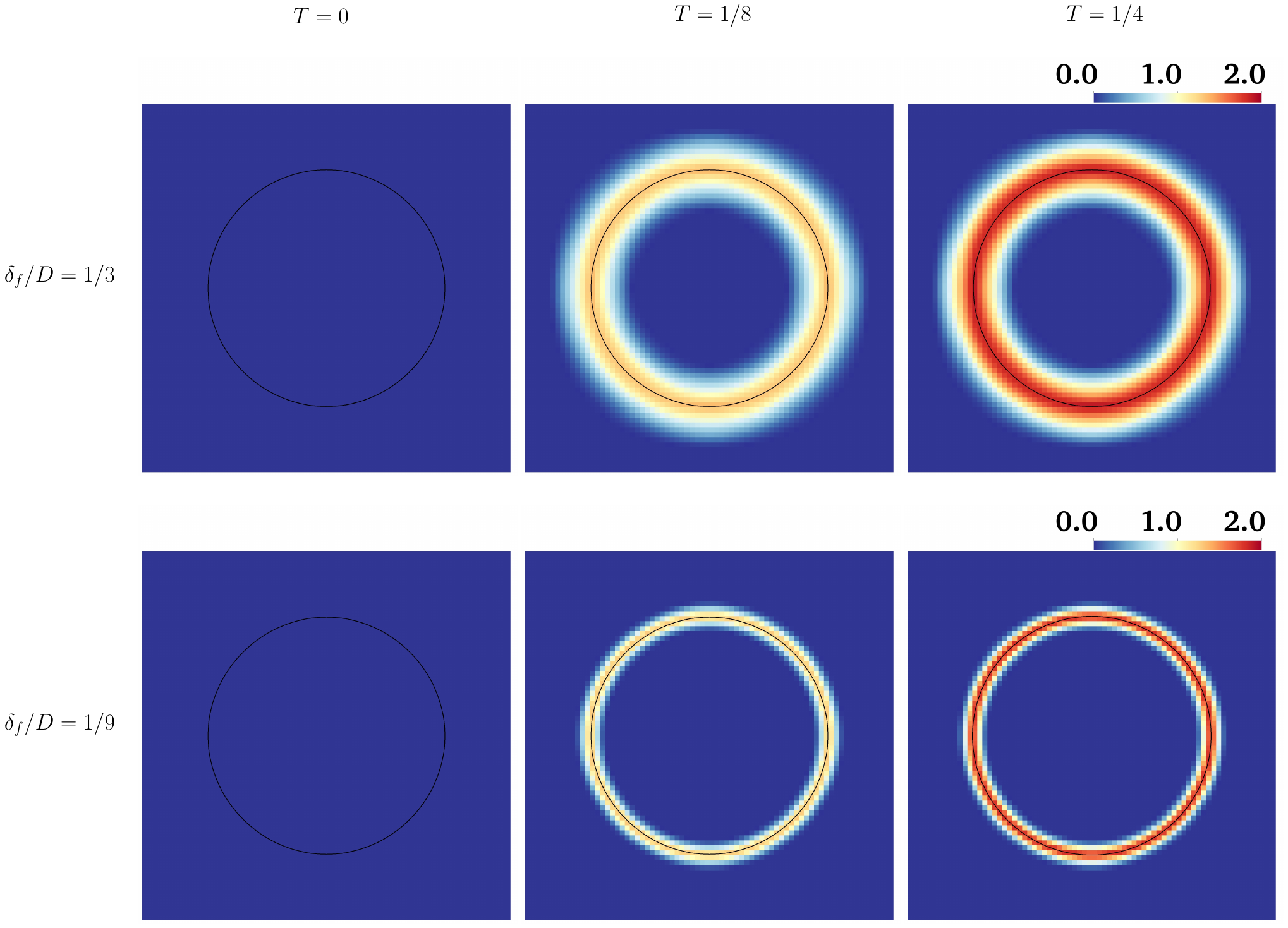}
	\caption{Isocontours of $\mid\delta_f F_I\mid$ at the absolute minimum value ($T = 0$) and at its absolute maximum value ($T = 1/4$). We also show $F_I$ at a middle point between the two extremes at $T = 1/8$. The forcing is sinusoidal in nature with respect to time. $F_I$ is shown at 2 different filter resolution of $\delta_f/D = 1/3$ and $\delta_f/D = 1/9$.}
	\label{fig:iso_fi}
\end{figure}

The filtered equation has four terms that govern the solution. The time derivative $\partial(\alpha\overline{u})/\partial t$ represents the rate of change of the solution in time. $\overline{\nabla G}\cdot\nabla(\alpha\overline{u})$ represents an advection term. The two other terms that govern the solution are the forcing term $F_I$ and the sub-filter scale term $\tau_\mathrm{sfs}$. We examine how $F_I$ varies with the filter width. The expression for $F_I$ is,
\begin{eqnarray}
	F_I = \overline{\nabla G}\cdot\iint_{\bm{y}\in S_I}u_I \bm{n}g(\bm{x}-\bm{y})dS = \overline{\nabla G}\cdot\left(u_I\iint_{\bm{y}\in S_I}\bm{n}g(\bm{x}-\bm{y})dS\right) = u_I\left(\overline{\nabla G}\cdot\nabla\alpha\right),\label{eq:FI_eqn}
\end{eqnarray}
where we have used the identity $\iint_{\bm{y}\in S_I}\bm{n}g(\bm{x}-\bm{y})dS = \nabla\alpha$ and the fact that $u_I$ does not vary along the boundary. Since $\nabla\alpha$ scales as $1/\delta_f$, equation (\ref{eq:FI_eqn}) shows that $F_I$ also scales with $1/\delta_f$. Figure \ref{fig:iso_fi} shows the isocontours of the forcing term $\mid\delta_f F_I\mid$ at varying $\delta_f/D$. The isocontours are shown when the forcing at the interface is at its maximum ($T = 1/4$) and minimum ($T = 0$) absolute magnitude as well as a time in between ($T = 1/8$). This term is computed in a purely analytical fashion using equation (\ref{eq:GE_3b}). In order to highlight how the filter width affects the forcing field, we show two different cases at $\delta_f/D = 1/3$ and $\delta_f/D = 1/9$. We observe that the forcing term $F_I$ is zero everywhere apart from a region within $\delta_f /2$ of the interface. This is inline with what would be expected since $F_I$ is computed by extrapolating the quantity at the interface onto the nearby cells within a region $\delta_f /2$ away from the interface. Furthermore, $F_I$ changes in a sinusoidal manner dependent on the boundary condition at the interface, dictated by $u_I$ which is defined by equation (\ref{eq:GE_3b}). The larger the filter width the more smeared $F_I$ is and hence has a larger area of influence to the field. As $\delta_f$ is reduced we see a sharper representation of $F_I$. While the interface is a sharp discontinuity, the volume-filtering approach removes the notion of the sharp boundary. Instead, the boundary condition is enforced by the forcing term $F_I$ over a narrow band whose width is dependent on the filter width $\delta_f$.

\begin{figure}
	\centering
	\includegraphics[width=0.9\linewidth]{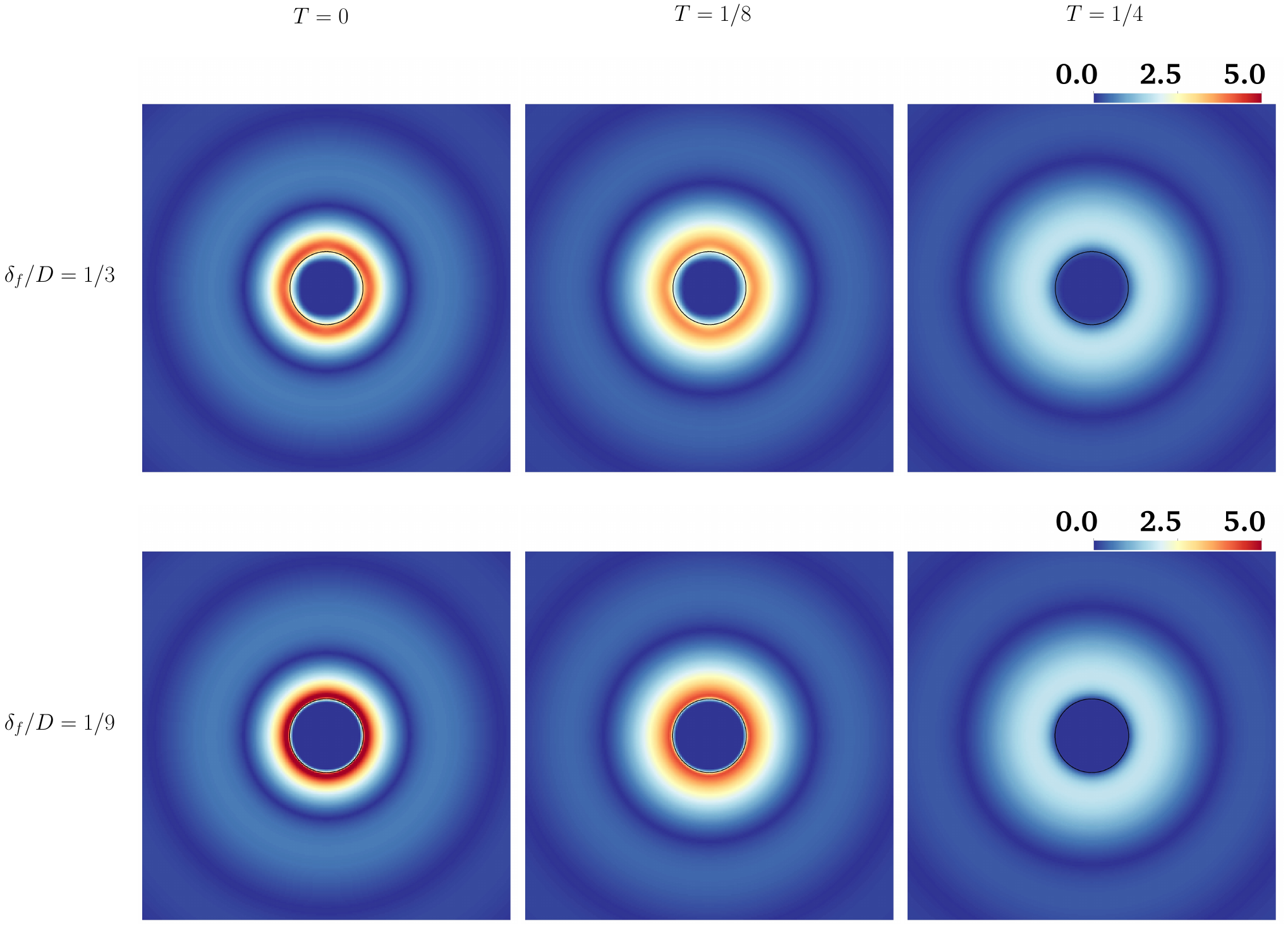}
	\caption{Isocontours of $\mid\tau_{\mathrm{sfs}}/\delta_f^2\mid$ is shown three different time snapshots. The results are computed purely analytically apart from the filtering procedure. We show two different filter widths, $\delta_f/D = 1/3$ and $\delta_f/D = 1/9$.}
	\label{fig:iso_sfs}
\end{figure}
$\tau_\mathrm{sfs}$ is the sub-grid term that comes out of the volume-filtering process and is given in equation ($\ref{eq:VF_4}$). Expanding $\tau_\mathrm{sfs}$ using a Taylor series gives us,
\begin{eqnarray}
	\tau_\mathrm{sfs} = \iiint_{\bm{y}\in\mathbbm{R}^3}\mathbbm{1}(\bm{y},t)\left(\nabla G(\bm{y})-\nabla G(\bm{x})\right)\cdot((\bm{y}-\bm{x})\cdot\nabla(\nabla u)(\bm{x}))g(\bm{x}-\bm{y})dV + \mathrm{HOT},\label{eq:sfs_taylor}
\end{eqnarray}
where $\mathrm{HOT}$ refers to the higher order terms in the Taylor series expansion. The equation above shows that $\tau_\mathrm{sfs}$ scales with $\delta_f^2$. Figure \ref{fig:iso_sfs} shows the isocontours for $\mid\tau_{\mathrm{sfs}}/\delta_f^2\mid$. The results are shown for both the extreme cases when $F_I$ is at its maximum amplitude and when $F_I$ is zero. The snapshots are at the same time as those shown in figure \ref{fig:iso_fi}. We observe that $\tau_{\mathrm{sfs}}$ is maximum near the interface and gradually reduces as we go away from it while following a wave-like pattern. This is because $\tau_{\mathrm{sfs}}$ is a function of $\overline{\nabla G}$ which is a function that decreases as we move away from the interface. The wave like pattern is due to the dot product between $\nabla G$ and $\nabla u$. At this resolution and interface sharpness, for both maximum and minimum absolute forcing at the interface the value of $\tau_\mathrm{sfs}$ is 3-4 orders of magnitude smaller than the advection term or the forcing term. From this, we infer that for a sharp enough interface the contribution of $\tau_\mathrm{sfs}$ is negligible and does not affect the accuracy of the filtered solution.

\begin{figure}\centering
  \begin{subfigure}{0.45\linewidth}
    \includegraphics[width=\linewidth]{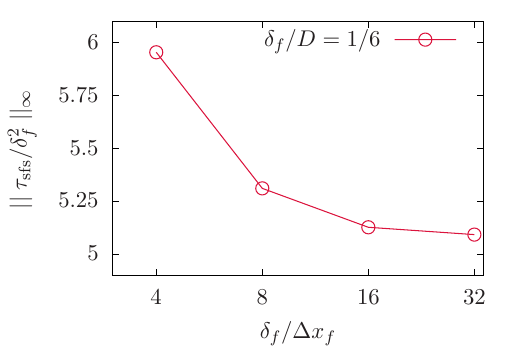}
    \caption{\label{fig:linf_dfdx_none}}
  \end{subfigure}
  \begin{subfigure}{0.45\linewidth}
    \includegraphics[width=\linewidth]{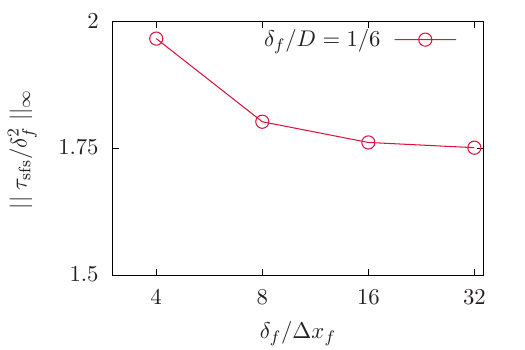}
    \caption{\label{fig:linf_dfdx_max}}
  \end{subfigure}
	\caption{$\mid\mid\tau_\mathrm{sfs}/\delta_f^2\mid\mid_\infty$ at $\delta_f/D = 1/6$ for varying $\delta_f/\Delta x_f$. We show the results are the two extreme time periods in the sinusoidal solution (a) $T = 1/2$ and (b) $T = 1/4$. We see that as we increase the number of subgrid points used for the filtering procedure, the value of $\mid\mid\tau_\mathrm{sfs}/\delta_f^2\mid\mid_\infty$ converges at $\delta_f/\Delta x_f = 32$ \label{fig:linf_dfdx}}
\end{figure}
Before comparing $\tau_\mathrm{sfs}$ to the other terms in the equation, we first observe how the $L_\infty$-norm of $\tau_\mathrm{sfs}$ varies with varying $\delta_f/\Delta x_f$. Here, $\Delta x_f$ is the subgrid mesh that we employ within the filtering procedure in order to remove numerical errors from the filtering procedure as much as possible and in order to compute an accurate volume fraction field and hence a more accurate $\tau_\mathrm{sfs}$ value at any point $\bm{x}$. We take a filter width of $\delta_f/D = 1/6$ and compute $\mid\mid\tau_\mathrm{sfs}/\delta_f^2\mid\mid_\infty$ and show the values at $T = 1/2$ and $T = 1/4$. Figure \ref{fig:linf_dfdx} shows that the value converges when $\delta_f/\Delta x_f = 32$.

\begin{figure}\centering
  \begin{subfigure}{0.45\linewidth}
    \includegraphics[width=\linewidth]{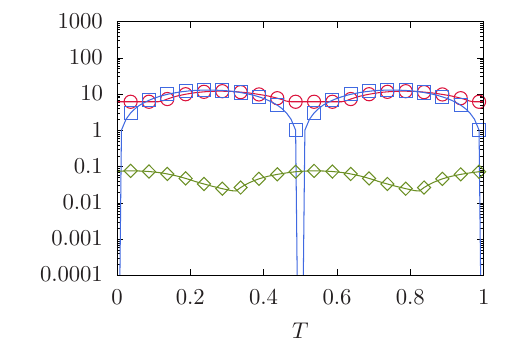}
    \caption{$\delta_f/D = 1/3$\label{fig:fig9a}}
  \end{subfigure}
  \begin{subfigure}{0.45\linewidth}
    \includegraphics[width=\linewidth]{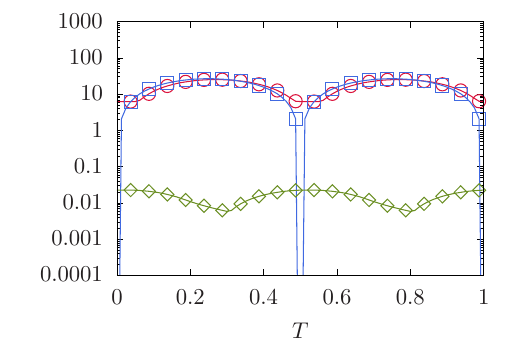}
    \caption{$\delta_f/D = 1/6$\label{fig:fig9b}}
  \end{subfigure}
	\begin{subfigure}{0.45\linewidth}
    \includegraphics[width=\linewidth]{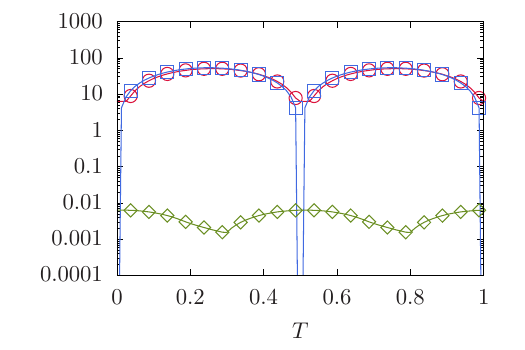}
    \caption{$\delta_f/D = 1/12$\label{fig:fig9c}}
  \end{subfigure}
	\begin{subfigure}{0.45\linewidth}
    \includegraphics[width=\linewidth]{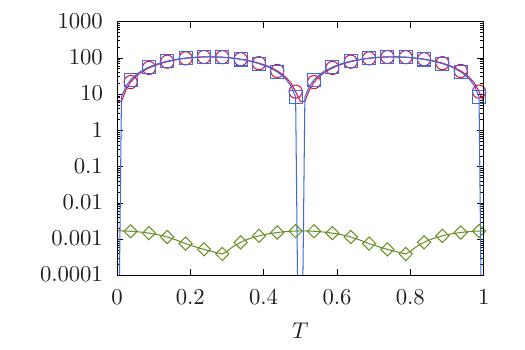}
    \caption{$\delta_f/D = 1/24$\label{fig:fig9d}}
  \end{subfigure}
	\caption{Time series of the $L_\infty$ norms is shown for the different terms in equation \ref{eq:VF_5}. $\mid\mid F_I\mid\mid_\infty$ denoted by (\squares). $\mid\mid\overline{\nabla G}\cdot\left(\alpha\overline{u}\right)\mid\mid_\infty$ is denoted by (\circles). Lastly $\mid\mid\tau_{\mathrm{sfs}}\mid\mid_\infty$ denoted by (\diamonds). The terms are computed using the analytical solution. We show the results for four different filter widths ranging from $\delta_f/D = 1/3$ to $\delta_f/D = 1/24$\label{fig:figure7}}
\end{figure}

To shed more light on how the contribution of $\tau_\mathrm{sfs}$ changes in comparison to the advection term and forcing term, figure \ref{fig:figure7} shows the $L_\infty$-norm of these terms different interface sharpness ranging from the coarsest interface sharpness of $(\delta_f/D)^{-1} = 3$ to the sharpest interface of $(\delta_f/D)^{-1} = 24$. We keep the subgrid mesh width constant at $\delta_f/\Delta x_f = 32$. The results clearly validate what is observed from the visual snapshots shown in figure \ref{fig:iso_sfs}. As we decrease $\delta_f$, $\tau_\mathrm{sfs}$ reduces, thus reducing its contribution to the solution. Furthermore, we see that $\tau_{\mathrm{sfs}}$ peaks when the forcing term $F_I$ is zero and is at its minimum when $F_I$ is at its maximum amplitude. Figure \ref{fig:figure7} shows that at the coarsest interface sharpness $\tau_{\mathrm{sfs}}$ is around 2 orders of magnitude smaller. As the interface sharpness is refined this number reduces to 4-5 orders of magnitude compared to the other two terms in the governing equations. We infer that for cases with coarser interface sharpness, a modeling approach would be warranted in order to obtain higher accuracy. However, for a fine enough interface sharpness, $\tau_{\mathrm{sfs}}$ is negligible.

\begin{figure}\centering
  \begin{subfigure}{0.45\linewidth}
    \includegraphics[width=\linewidth]{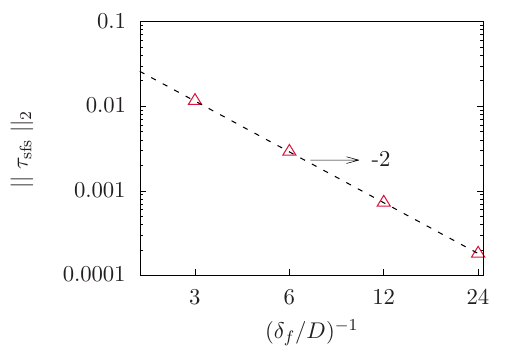}
    \caption{\label{fig:fig10a}}
  \end{subfigure}
  \begin{subfigure}{0.45\linewidth}
    \includegraphics[width=\linewidth]{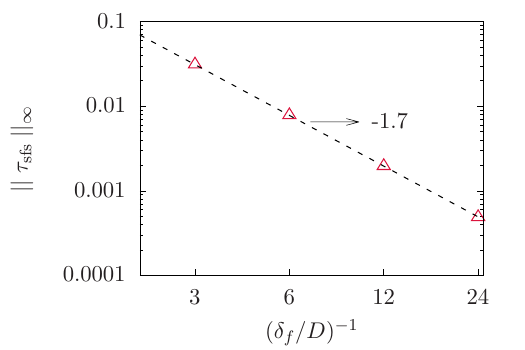}
    \caption{\label{fig:fig10b}}
  \end{subfigure}
	\begin{subfigure}{0.45\linewidth}
    \includegraphics[width=\linewidth]{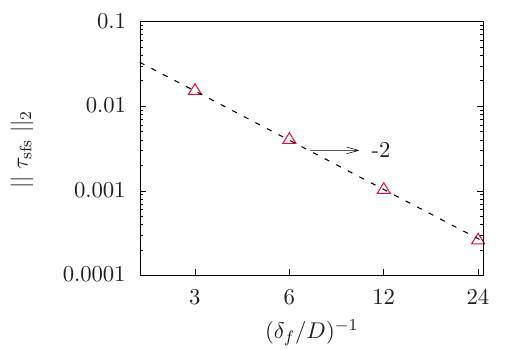}
    \caption{\label{fig:fig10c}}
  \end{subfigure}
	\begin{subfigure}{0.45\linewidth}
    \includegraphics[width=\linewidth]{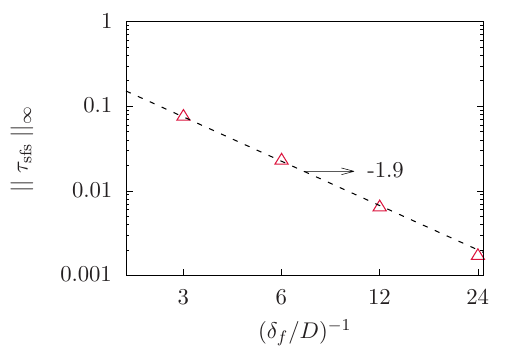}
    \caption{\label{fig:fig10d}}
  \end{subfigure}
	\caption{$\mid\mid\tau_\mathrm{sfs}\mid\mid_2$ and $\mid\mid\tau_\mathrm{sfs}\mid\mid_\infty$ at $T = 1/4$ (top) and $T = 1/2$ (bottom). The results are shown for varying filter width ($\delta_f/D$) while keeping $\delta_f/\Delta x_f = 32$. The terms are computed analytically.\label{fig:figure10}}
\end{figure}

In order to confirm that $\tau_\mathrm{sfs}$ scales with $\delta_f^2$, we plot $\mid\mid\tau_\mathrm{sfs}\mid\mid_2$ and $\mid\mid\tau_\mathrm{sfs}\mid\mid_\infty$ for increasing $(\delta_f/D)^{-1}$, i.e. increasing sharpness. In doing so, we calculate the order of convergence with respect to $(\delta_f/D)^{-1}$. The results are shown when $\tau_\mathrm{sfs}$ is at its absolute maximum at $T = 1/2$ and when it is at its absolute minimum at $T = 1/4$ in figure \ref{fig:figure10}. We show the order of convergence for the $\mid\mid\tau_\mathrm{sfs}\mid\mid_2$ and $\mid\mid\tau_\mathrm{sfs}\mid\mid_\infty$ norm when $\tau_\mathrm{sfs}$ is at its maximum and when $\tau_\mathrm{sfs}$ is at its minimum. Similar to above the ratio of $\delta_f/\Delta x_f = 32$ is kept constant. We observe that $\tau_{\mathrm{sfs}}$ scales with $\delta_f^2$, hence giving us a second order rate of convergence for $\tau_{\mathrm{sfs}}$ with $\delta_f$.

\subsection{Aposteriori analysis}
\label{sssec:apos_analysis}

\begin{figure}
	\centering
	\includegraphics[width=0.9\linewidth]{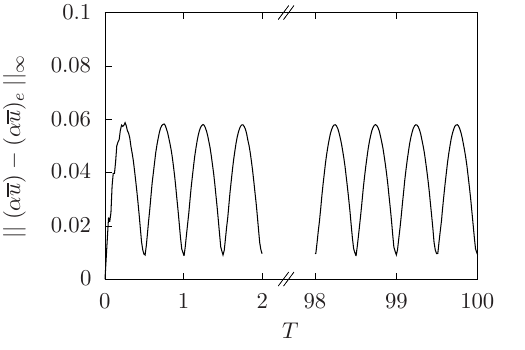}
	\caption{$\mid\mid(\alpha\overline{u})-(\alpha\overline{u})_e\mid\mid_\infty$ for 100 time periods in order to show stability of the numerical solution. The result is run at $\delta_f/D = 1/12$, $\delta_f/\Delta x = 4$ and a $\mathrm{CFL} = 0.5$. The subgrid mesh is kept at $\delta_f/\Delta x_f = 32$}
	\label{fig:stability}
\end{figure}
We now solve numerically the volume-filtered equations derived in section \ref{sssec:derivation_vfib}, and compare the numerical solution with the filtered analytical solution $(\alpha\overline{u})_e$. For now, we keep $\tau_\mathrm{sfs}$ turned off.

In order to show the stability of the method we run the case at a filter width of $\delta_f/D = 1/12$ and a spatial resolution of $\delta_f/\Delta x = 4$. The simulation is run at a Courants-Friedrichs Lewy (CFL) number of $CFL = 0.5$. We run the simulation for a total of 100 periods and we report the $\mid\mid(\alpha\overline{u})-(\alpha\overline{u})_e\mid\mid_\infty$ error in figure \ref{fig:stability}. The results show a constant sinusoidal graph, showing that the simulation is stable for long periods of time and does not blow up numerically. Furthermore, we observe that the error is at its maximum when the forcing term is greatest, at points $T = 1/4$ and $T = 3/4$. The error is at its minimum when the forcing is zero at the interface, $T = 0.5$ and $T = 1$.
\begin{figure}
	\centering
	\includegraphics[width=0.9\linewidth]{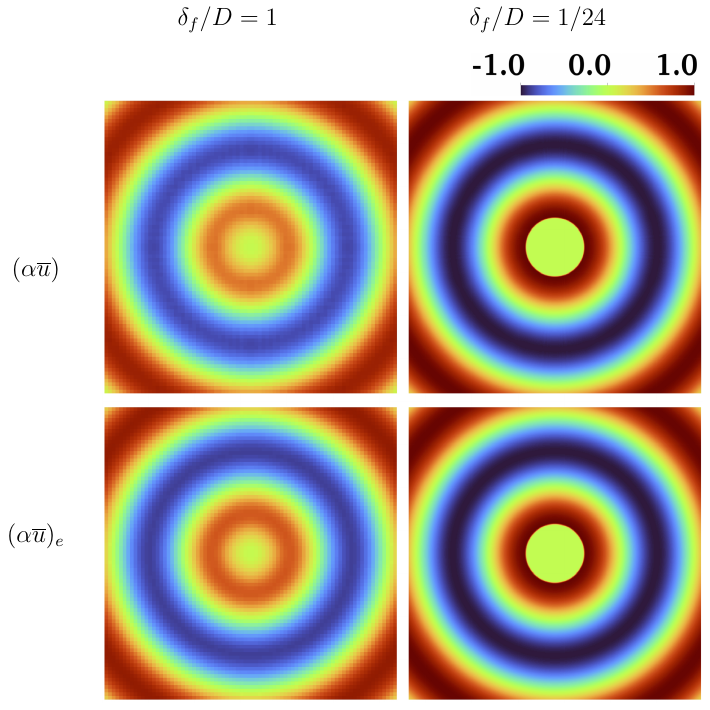}
	\caption{Isocontours of the numerical solution $(\alpha\overline{u})$ and the analytically filtered solution $(\alpha\overline{u})_e$ for filter width $\delta_f/D = 1$ (left) and $\delta_f/D = 1/24$ (right). The results are shown for phase $T = 1/4$. The simulations are run at $CFL = 0.25$ and the subgrid mesh is kept such that $\delta_f/\Delta x_f = 32$. The main grid spatial resolution is set at $\delta_f/\Delta x = 16$.}
	\label{fig:num_iso}
\end{figure}

\begin{figure}\centering
  \begin{subfigure}{0.45\linewidth}
    \includegraphics[width=\linewidth]{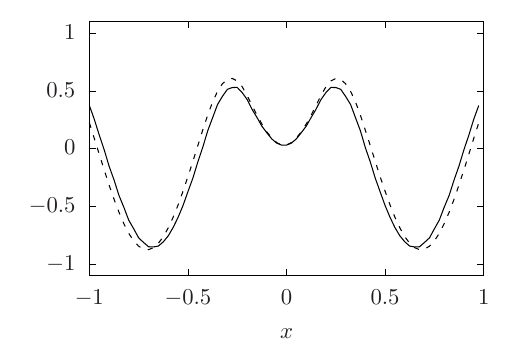}
    \caption{$\delta_f/D = 1$\label{fig:fig14a}}
  \end{subfigure}
  \begin{subfigure}{0.45\linewidth}
    \includegraphics[width=\linewidth]{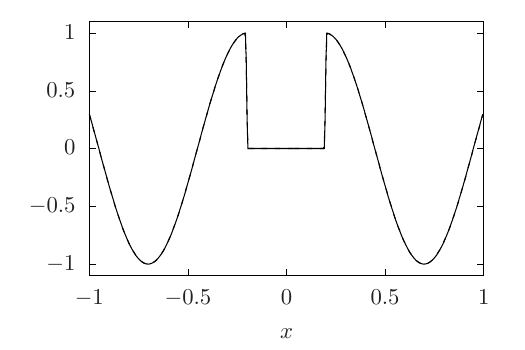}
    \caption{$\delta_f/D = 1/24$\label{fig:fig14b}}
  \end{subfigure}
	\caption{Value of $(\alpha\overline{u})$ (\stl) vs $x$ at $T = 1/4$. For comparison we also plot the analytically filtered solution $(\alpha\overline{u})_e$ (\dedash). The result shown is a line cut horizontally and centered in the vertical axis. We show the results for a coarse filter width of $\delta_f/D = 1$ and a fine filter width of $\delta_f/D = 1/24$. \label{fig:figure14}}
\end{figure}

Figure \ref{fig:num_iso} shows the isocontours of $(\alpha\overline{u})$ and $(\alpha\overline{u})_e$ at $\delta_f/D = 1$ and $\delta_f/D = 1/24$. We show the results at $T = 1/4$. both cases are run at a $CFL = 0.25$, $\delta_f/\Delta x = 4$. We make sure to maintain the subgrid mesh at $\delta_f/\Delta x_f = 32$. The numerical solution approaches the filtered analytical solution as we decrease the filter width $\delta_f/D$. For a more quantitative comparison, figure \ref{fig:figure14} plots the value of $(\alpha\overline{u})$ and $(\alpha\overline{u})_e$ as a line cut horizontally and centered in the vertical axis. We show the results for both a coarse filter width ($\delta_f/D = 1$) and a fine filter width ($\delta_f/D = 1/24$). We clearly see that as the filter width is reduced, the numerical solution is much closer to the filtered analytical solution. This is because as we coarsen the filter, the contribution from the sub-filter scale term, $\tau_\mathrm{sfs}$ increases. Hence, for larger filter widths, Errors can be expected unless a modeling approach is utilized in order to account for the contributuon from $\tau_\mathrm{sfs}$. For finer filter widths, the contribution of $\tau_\mathrm{sfs}$ is negligible, and therefore neglecting the term does not affect the quality of the solution.

\begin{figure}\centering
  \begin{subfigure}{0.45\linewidth}
    \includegraphics[width=\linewidth]{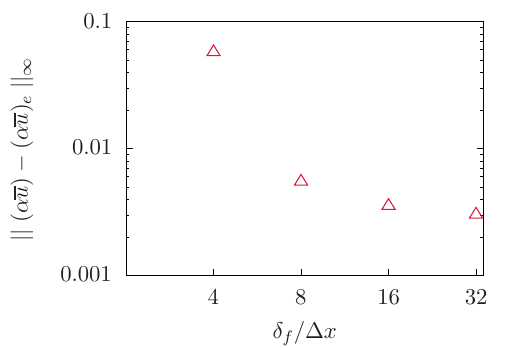}
    \caption{\label{fig:max_converge}}
  \end{subfigure}
  \begin{subfigure}{0.45\linewidth}
    \includegraphics[width=\linewidth]{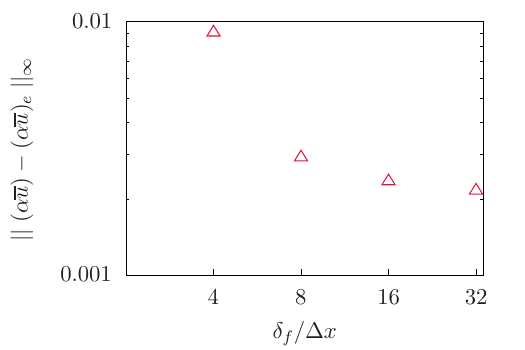}
    \caption{\label{fig:none_converge}}
  \end{subfigure}
	\caption{$\mid\mid(\alpha\overline{u})-(\alpha\overline{u})_e\mid\mid_\infty$ norm for varying $\delta_f/\Delta x$ at $(\delta_f/D)^-1 = 1/12$. (a) shows the error when the interface forcing is at its maximum and (b) show the error when the interface forcing is at its minumum. We show that the error converges to stationary value as we increase the number of grid points. \label{fig:figure_converge}}
\end{figure}

\begin{figure}\centering
  \begin{subfigure}{0.45\linewidth}
    \includegraphics[width=\linewidth]{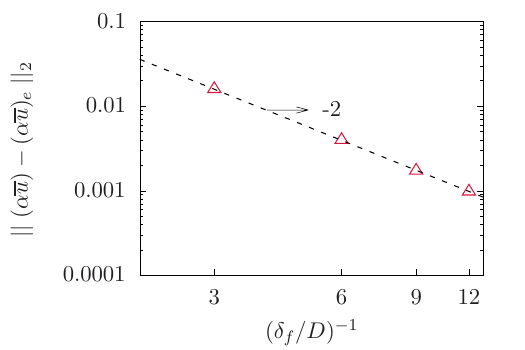}
    \caption{\label{fig:fig13a}}
  \end{subfigure}
  \begin{subfigure}{0.45\linewidth}
    \includegraphics[width=\linewidth]{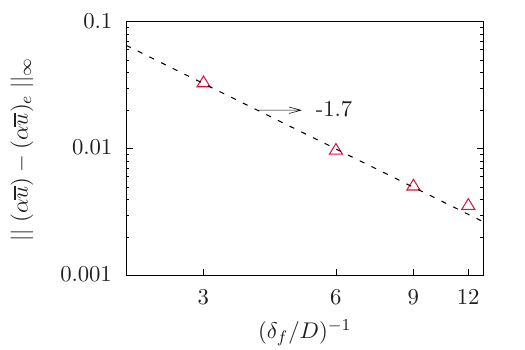}
    \caption{\label{fig:fig13b}}
  \end{subfigure}
	\begin{subfigure}{0.45\linewidth}
    \includegraphics[width=\linewidth]{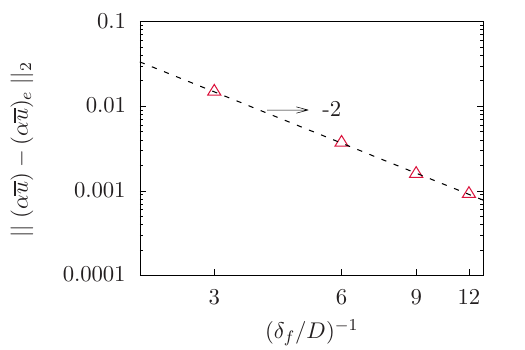}
    \caption{\label{fig:fig13c}}
  \end{subfigure}
	\begin{subfigure}{0.45\linewidth}
    \includegraphics[width=\linewidth]{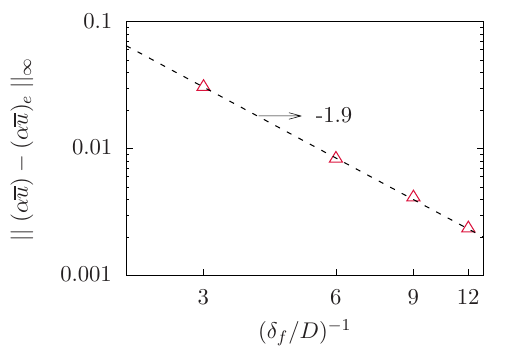}
    \caption{\label{fig:fig13d}}
  \end{subfigure}
	\caption{$\mid\mid(\alpha\overline{u})-(\alpha\overline{u})_e\mid\mid_2$ and $\mid\mid(\alpha\overline{u})-(\alpha\overline{u})_e\mid\mid_\infty$ norm when the interface forcing $F_I$ is at its absolute maximum (top) and when it is at its absolute minimum (bottom). The results are shown for increasing $(\delta_f/D)^{-1}$ while keeping $\delta_f/\Delta x = 16$. All simulations are run at $\mathrm{CFL} = 0.1.$\label{fig:figure12}}
\end{figure}

In order to accurately calculate the order of convergence for $(\alpha\overline{u})$, we first need to eliminate all numerical errors as much as possible. One of the spatial parameters that govern the numerical solution is the ratio between the filter with and the mesh spacing $\delta_f/\Delta x$. The subgrid filtering mesh is kept constant at $\delta_f/\Delta x_f = 32$ and $\tau_\mathrm{sfs}$ is still turned off. Figure \ref{fig:figure_converge} shows $\mid\mid(\alpha\overline{u})-(\alpha\overline{u})_e\mid\mid_\infty$ for varying $\delta_f/\Delta x$ ranging from $\delta_f/\Delta x = 4$ to $\delta_f/\Delta x = 32$. We see from the graph that the results converge at $\delta_f/\Delta x = 16$. We obtain the order of convergence for varying $\delta_f/D$. In order to highlight how interface sharpness affects the results, figure \ref{fig:figure12} shows the $\mid\mid(\alpha\overline{u})-(\alpha\overline{u})_e\mid\mid_2$ and $\mid\mid(\alpha\overline{u})-(\alpha\overline{u})_e\mid\mid_\infty$ error. All simulations are run at $\mathrm{CFL} = 0.1$. The ratio of the filter width to the mesh width is fixed at $\delta_f/\Delta x = 16$. Currently we turn off $\tau_\mathrm{sfs}$ and only consider the advection and forcing term. Looking at the order of convergence for the $L_2$ and $L_\infty$ norm error we observe that the slope obtained is similar to that of the order of convergence for $\tau_\mathrm{sfs}$. This shows that the VF-IB method is a stable second order method and performs well even under spatially coarse conditions.

\begin{figure}
	\centering
	\includegraphics[width=\linewidth]{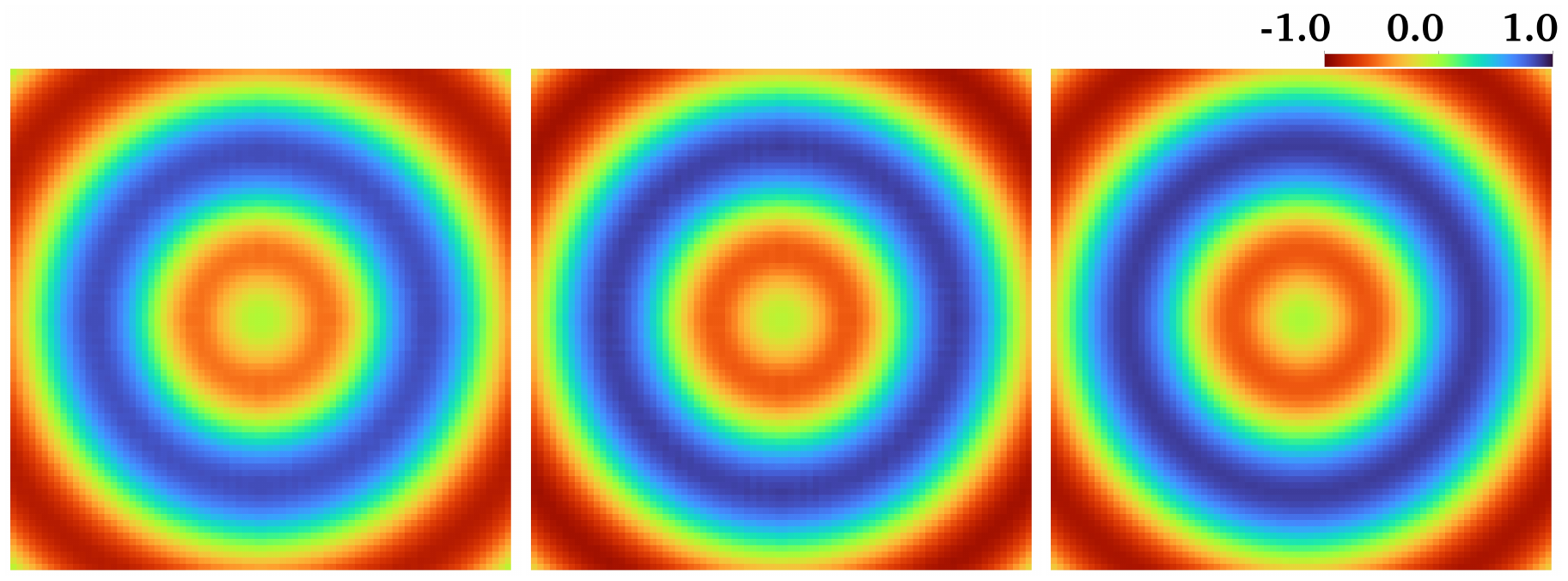}
	\caption{Isocontours of $(\alpha\overline{u})$ when $\tau_\mathrm{sfs}$ is turned off (left) and then $\tau_\mathrm{sfs}$ is turned on (middle). We also show the filtered analytical solution $(\alpha\overline{u})_e$ for comparison (right). The results are shown at the filter width of $\delta_f/D = 1$ to show how including $\tau_\mathrm{sfs}$ makes the numerical solution closer to the filtered analytical solution at coarse filter sizes.}
	\label{fig:iso_contour_sfs}
\end{figure}

\begin{figure}\centering
  \begin{subfigure}{0.45\linewidth}
    \includegraphics[width=\linewidth]{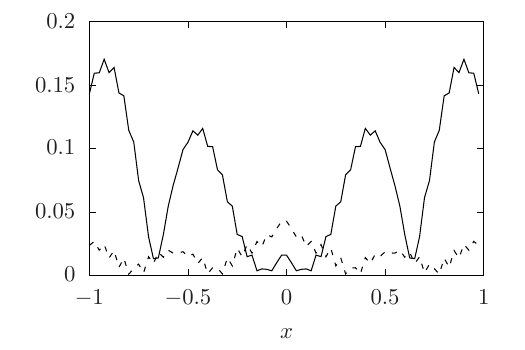}
    \caption{\label{fig:sfs_quant_max}}
  \end{subfigure}
  \begin{subfigure}{0.45\linewidth}
    \includegraphics[width=\linewidth]{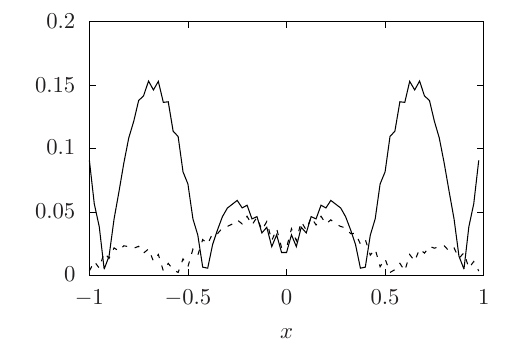}
    \caption{\label{fig:sfs_quant_min}}
  \end{subfigure}
	\caption{Value of $\mid\mid(\alpha\overline{u}) - (\alpha\overline{u})_e\mid\mid_\infty$ vs $x$ at (a) $T = 1/4$ and (b) $T = 1/2$. The result shown is a line cut horizontally and centered in the vertical axis. We show the result when $\tau_\mathrm{sfs}$ is turned off (\stl) and when $\tau_\mathrm{sfs}$ is turned on (\dd). We run the simulations at $\delta_f/D = 1$. \label{fig:sfs_quant}}
\end{figure}

The contribution of the sub-filter scale term $\tau_\mathrm{sfs}$ on the numerical solution depends on $\delta_f/D$. For large $\delta_f/D$ values, the contribution from $\tau_\mathrm{sfs}$ is large enough that including it significantly improves the accuracy of the solution compared to the filtered analytical solution $(\alpha\overline{u})_e$. Figure \ref{fig:iso_contour_sfs} shows the isocontours of $(\alpha\overline{u})$ when $\tau_\mathrm{sfs}$ is turned on and when we ignore it. We also show the filtered analytical solution $(\alpha\overline{u})_e$ for comparison. The results are shown at a coarse filter width of $\delta_f/D = 1$. The subgrid mesh is kept constant at $\delta_f/\Delta x_f = 32$ and the main grid resolution is $\delta_f/\Delta x = 16$. The cases are run at $\mathrm{CFL} = 0.1$. From the results, we observe that when we include $\tau_\mathrm{sfs}$ the error is significantly lower than when it is excluded. This shows that for coarser filter widths, including $\tau_\mathrm{sfs}$ helpes improve the accuracy of the solution. In order to show this in a quantitative manner, figure \ref{fig:sfs_quant} shows the value of $\mid\mid(\alpha\overline{u}) - (\alpha\overline{u})_e\mid\mid_\infty$ as a line cut horizontally and centered in the vertical direction. We show the results, both when $\tau_\mathrm{sfs}$ is turned off and then it is turned on. This is further proof that inclusion of $\tau_\mathrm{sfs}$ at coarser filter widths helps produce a more accurate numerical solution.

\begin{figure}\centering
  \begin{subfigure}{0.45\linewidth}
    \includegraphics[width=\linewidth]{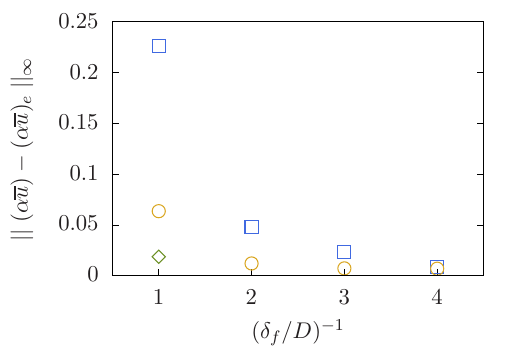}
    \caption{\label{fig:max_sfs}}
  \end{subfigure}
  \begin{subfigure}{0.45\linewidth}
    \includegraphics[width=\linewidth]{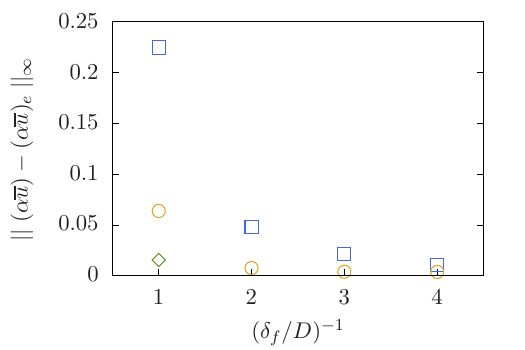}
    \caption{\label{fig:mim_sfs}}
  \end{subfigure}
	\caption{$\mid\mid(\alpha\overline{u})-(\alpha\overline{u})_e\mid\mid_\infty$ norm for varying $(\delta_f/D)^{-1}$ at $\delta_f/\Delta x = 16$. (\squares) correspond to the result when $\tau_\mathrm{sfs}$ is turned off and (\circles) correspond to the result when $\tau_\mathrm{sfs}$ is turned on. In order to show how $D/\Delta x$ affects the solution, (\diamonds) shows the results at $D/\Delta x = 32$ when $\tau_\mathrm{sfs}$ is turned on.	(a) shows the error when the interface forcing is at its maximum and (b) show the error when the interface forcing is at its minumum. We show that for coarse filter widths, turning on $\tau_\mathrm{sfs}$ helps significantly reduce the error. This error can be further reduced by resolving the interface grid resolution $D/\Delta x$. \label{fig:figure_sfson}}
\end{figure}

Figure \ref{fig:figure_sfson} shows how $\mid\mid(\alpha\overline{u})-(\alpha\overline{u})_e\mid\mid_\infty$ varies with changing $\delta_f/D$. The results are run at $\delta_f/\Delta x = 16$ and a subgrid resolution of $\delta_f/\Delta x_f = 32$. In order to highlight how refining the interface grid resolution $D/\Delta x$ affects the solution we run an extra case at $D/\Delta x = 32$ at $\delta_f/D = 1$ when $\tau_\mathrm{sfs}$ is turned on. The results show that as we reduce $\delta_f/D$, the contribution from $\tau_\mathrm{sfs}$ significantly reduces. We observe that at $\delta_f/D = 1/4$ the error is negligible enough. At this resolution or finer resolutions, $\tau_\mathrm{sfs}$ can be ignored without any reduction in the accuracy of the solution. Furthermore, we show the importance of resolving the interface grid resolution. At a $\delta_f/D = 1$ we witness that the error is significantly lower at $D/\Delta x = 32$ than it is at $D/\Delta x = 16$. We therefore conclude that for coarser filter resolutions modeling $\tau_\mathrm{sfs}$ will improve the accuracy, however as we refine $\delta_f/D$, we can ignore the term. Furthermore, refining the grid resolution $D/\Delta x$ will also improve the accuracy of the solution.

\section{Conclusion}
\label{sec:sec_conclusion}

In this paper we extend the novel Volume-Filtered Immersed Boundary method (VF-IB) by \citet{daveVolumefilteringImmersedBoundary2023} through characterizing the unclosed term $\tau_\mathrm{sfs}$ and the forcing term $F_I$ for varying spatial parameters. This helps us better understand how unclosed terms play a role in the VF-IB method and how varying filter width affects the quality of the solution.

In order to show this, we take the case of a varying coefficient hyperbolic equation. In a purely theoretical fashion, through the process of volume filtering \citep{andersonFluidMechanicalDescription1967} we are able to convert the point-wise equations into filtered transport equations. The immersed boundary plays a role in the forcing term, which appears as a surface integral on the right hand side of the transport equation. We show that using the volume-filtering process within an immersed boundary framework, we are able to theoretically transform the PDE that involves boundary conditions on surfaces into a filtered equation to solve, where the surface integrals define the boundary conditions at the interface. This process is not tied to any discretization and is viable for any topologically complex boundary. Furthermore, through the process of volume-filtering arise the sub-filter scale terms, namely $\tau_\mathrm{sfs}$. This paper studies the effect that interface sharpness has on the magnitude of $\tau_\mathrm{sfs}$, in comparison to the terms that are present in the volume-filtered transport equations. Furthermore, we also explore how $\tau_\mathrm{sfs}$ contributes to the solution for varying $\delta_f/D$.

We first study the seperate terms shown in the filtered transport equations in a purely analytical fashion as shown in section \ref{sssec:ap_analysis}. We see that for finer interface sharpness values the contribution of $\tau_\mathrm{sfs}$ is small, and can be considered to be negligible. Hence,  the term can be ignored without any significant reduction in accurary of the solution. As the interface sharpness becomes coarse this contribution increases, making $\tau_\mathrm{sfs}$ significant. At these coarse resolutions a modeling approach may be warranted in order to output high orders of accuracy. Furthermore, we take a look at the order of convergence for $\tau_\mathrm{sfs}$ and show that the sub-filter scale term scales with $\delta_f^2$, making the VF-IB method theoretically second order. We then conduct an aposteriori analysis as shown in section \ref{sssec:apos_analysis} and perform the simulations numerically using the VF-IB method. We initially ignore the effects of $\tau_\mathrm{sfs}$ and turn off the term. Again, we conduct simulations from the coarsest interface sharpness of $\delta_f/D = 1/3$ and gradually increase the sharpness to $\delta_f/D = 1/12$. The results are compared to the analytically filtered solution and the results are as expected. We observe a second order rate of convergence that is in line with what was observed during the apriori analysis. For sharp interface cases such as at $\delta_f/D = 1/12$. The results are comparable to the second order results in the cut-cell method by \citet{bradyFoundationsHighorderConservative2021}. Hence, we can say that for a sufficiently sharp interface, $\tau_\mathrm{sfs}$ can be ignored while still keeping obtaining a high order of accuracy.

In order to obtain high level of accuracy at coarser interface sharpness, the sub-filter scale terms would need to be modeled. To show this we take cases with varying $\delta_f/D$ from $\delta_f/D = 1$ to $\delta_f/D = 1/4$ and include $\tau_\mathrm{sfs}$ as part of the solution. From this we conclude that at extremely coarse resolution, $\tau_\mathrm{sfs}$ has a contribution large enough that it cannot be ignored and would have to be modeled in order to achieve good accuracy. However, as we refine $\delta_f/D$ this value quickly reduces in magnitude and does not play a major role in the solution. Keep in mind, these resolutions are coarse enough that they  can run at an extremely cheap computational cost and obtaining sub-filter scale models can help achieve robustness of the method while still obtaining rather good accuracy against the filtered analytical solution. Furthermore, we also shed light on the interface grid resolution $D/\Delta x$. We show that increasing $D/\Delta x$ significantly improves the results. A combination of accurate $\tau_\mathrm{sfs}$ models together with reasonable grid resolution can provide a robust method that is accurate while keeping the computational cost low.

Finally, we show that this methodical volume-filtering process can be conducted on any arbitrary hyperbolic PDE that involves boundary conditions on surfaces. The process allows us to solve for the transport equations on a cartesian grid using lagrangian markers as a tool to contruct the interface. This allows us to solve around any topologically complex interface that can be moving or static.

\section*{Acknowledgement}
The authors acknowledge support from the US National Science Foundation (award \#2028617, CBET-FD and EEC-GOALI award \#2216969). Computing resources were provided by Research Computing at Arizona State University.

\appendix

\section{Convergence analysis for embedded circle placed at different locations}\label{sec:appendix1}

\begin{figure}
	\centering
	\includegraphics[width=0.9\linewidth]{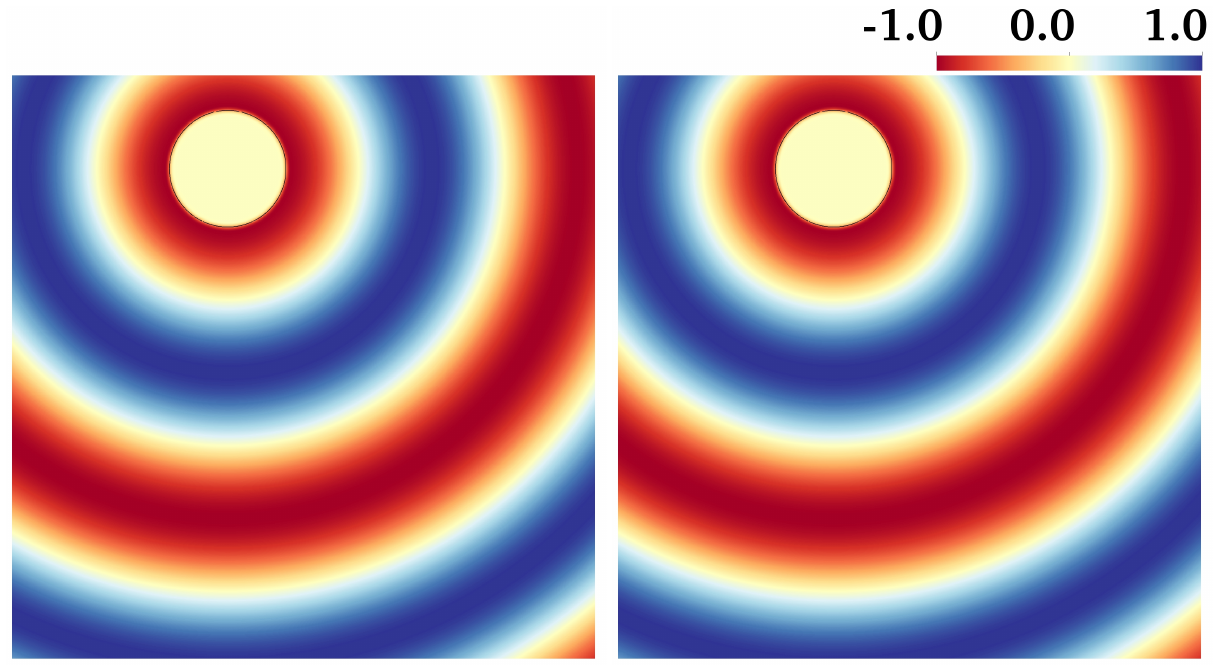}
	\caption{Isocontours of the numerical solution $(\alpha\overline{u})$ and the analytically filtered solution $(\alpha\overline{u})_e$ for filter width $\delta_f/D = 1/12$. The results are shown for phase $T = 3/4$. The simulations are run at $CFL = 0.1$ and the subgrid mesh is kept such that $\delta_f/\Delta x_f = 32$. The main grid spatial resolution is set at $\delta_f/\Delta x = 16$.}
	\label{fig:appendix_contour}
\end{figure}

\begin{figure}\centering
  \begin{subfigure}{0.45\linewidth}
    \includegraphics[width=\linewidth]{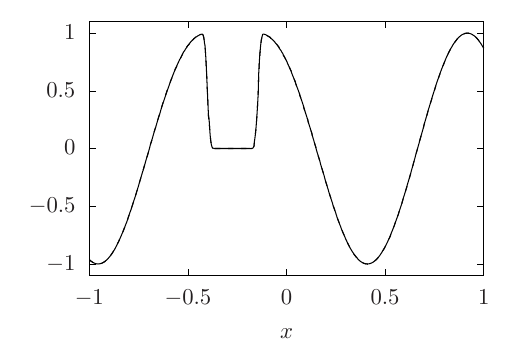}
    \caption{\label{fig:figaa}}
  \end{subfigure}
  \begin{subfigure}{0.45\linewidth}
    \includegraphics[width=\linewidth]{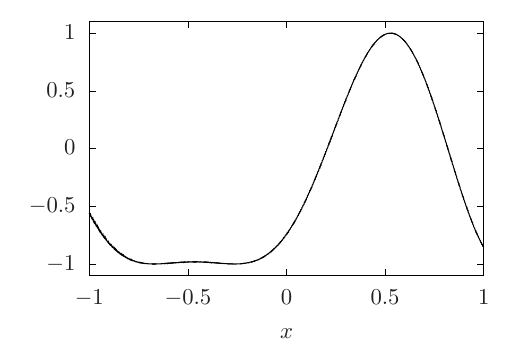}
    \caption{\label{fig:figab}}
  \end{subfigure}
	\caption{Value of $(\alpha\overline{u})$ (\stl) vs $x$ at $T = 3/4$. For comparison we also plot the analytically filtered solution $(\alpha\overline{u})_e$ (\dedash). The results shown are two lines picked at random that cut accross the domain. The first line goes from points $(-1,0.6)\rightarrow (1,0.4)$ (left) and the second line goes from points $(-1,0.25)\rightarrow (1,-0.53)$ (right). We show the results for a filter resolution of $\delta_f/D = 1/12$. \label{fig:plot_appendix}}
\end{figure}

Here we show the results when the circle is placed at 4 different locations. The center locations are such that, $l_1 = \left(0.5, 0.5\right)$ is within the top left quadrant. The other three locations are picked at randon such that, $l_2 = \left(0.5, -0.34\right)$, $l_3 = \left(-0.61, -0.43\right)$ and $l_4 = \left(-0.26, 0.68\right)$. Figure \ref{fig:appendix_contour} shows the isocontours of $(\alpha\overline{u})$ in comparison to $(\alpha\overline{u})_e$ when the circle is placed at location $l_4$. The simulation is run at a filter resolution of $\delta_f/D = 1/12$ and a spatial resolution of $\delta_f/\Delta x = 16$. We can visually compare the results and observe that the change in position of the circle does not affect the accuracy of the simulation. To further show this, figure \ref{fig:plot_appendix} plots the values of $(\alpha\overline{u})$ and $(\alpha\overline{u})_e$ by `cutting' accross the domain using two randomly selected lines. The results show that the numerical solution matches the analytically filtered solution well, thus making the method independent of the immersed boundary placement.

\begin{figure}\centering
  \begin{subfigure}{0.45\linewidth}
    \includegraphics[width=\linewidth]{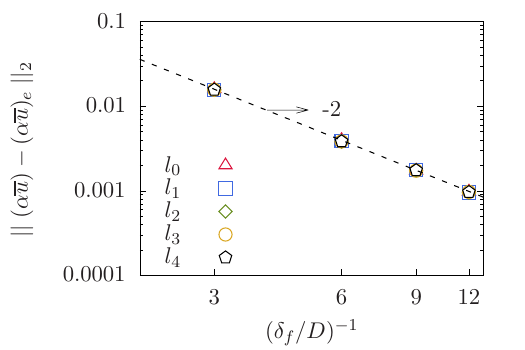}
    \caption{\label{fig:fig13a}}
  \end{subfigure}
  \begin{subfigure}{0.45\linewidth}
    \includegraphics[width=\linewidth]{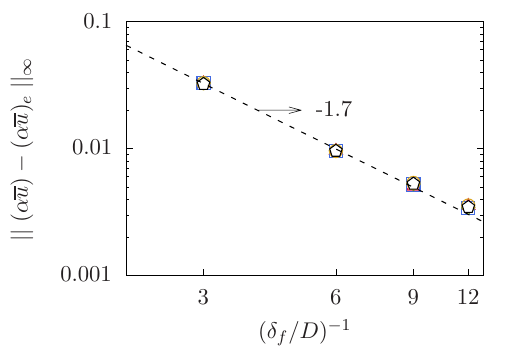}
    \caption{\label{fig:fig13b}}
  \end{subfigure}
	\begin{subfigure}{0.45\linewidth}
    \includegraphics[width=\linewidth]{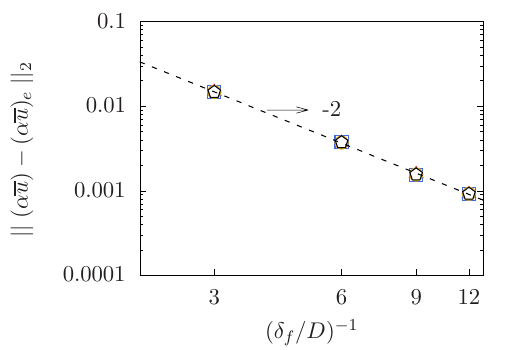}
    \caption{\label{fig:fig13c}}
  \end{subfigure}
	\begin{subfigure}{0.45\linewidth}
    \includegraphics[width=\linewidth]{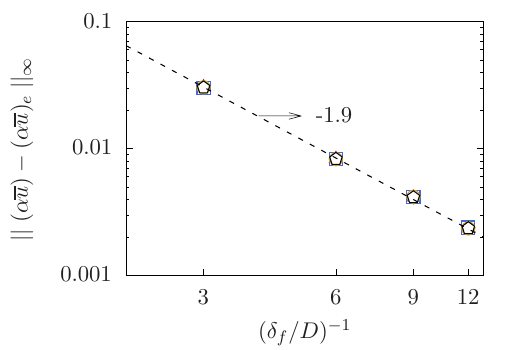}
    \caption{\label{fig:fig13d}}
  \end{subfigure}
	\caption{$\mid\mid(\alpha\overline{u})-(\alpha\overline{u})_e\mid\mid_2$ and $\mid\mid(\alpha\overline{u})-(\alpha\overline{u})_e\mid\mid_\infty$ norm when the interface forcing $F_I$ is at its absolute maximum (top) and when it is at its absolute minimum (bottom). The results are shown for increasing $(\delta_f/D)^{-1}$ while keeping $\delta_f/\Delta x = 16$. All simulations are run at $\mathrm{CFL} = 0.1$. We show the results for the circle placed at 4 different locations in comparison to the original position.\label{fig:fig_appendix}}
\end{figure}

Finally we perform a convergence analysis. For comparison we also show the results when the circle is in the center of the domain, $l_0 = \left(0,0\right)$. Figure \ref{fig:fig_appendix} shows that we get a second order convergence and the results are similar even when the IB is placed in different locations of the domain. Hence confirming that the second order convergence is not dependent on the location of the immersed boundary.

\bibliography{references_himanshu}
\end{document}